\documentclass{tac}
\usepackage{amsmath,amssymb,mathrsfs}
\usepackage{tikz-cd}
\usepackage{quiver}
\usepackage[noabbrev,capitalise]{cleveref}
\crefformat{equation}{(#2#1#3)}
\crefformat{enumi}{(#2#1#3)}
\crefformat{enumii}{(#2#1#3)}
\usepackage[
backend=biber,
sorting=nyt,
style=ieee
]{biblatex}
\usepackage{enumitem}
\setlist{listparindent = \parindent, parsep=0pt,}
\setenumerate[1]{label = (\roman*), ref = (\roman*)}
\setenumerate[2]{label = (\alph*), ref = (\alph*)}
\usepackage[X2,T1]{fontenc} 
\DeclareTextSymbolDefault{\CYRABHDZE}{X2}
\DeclareTextSymbolDefault{\cyrabhdze}{X2}
\newcommand{\Ezh}{\CYRABHDZE}

\newcommand{\GC}{\text{\rm\Ezh}}

\usepackage[colorlinks=true]{hyperref}
\hypersetup{allcolors=[rgb]{0.1,0.1,0.4}}

\author{Joshua L.\ Wrigley}

\address{Université Paris Cité, CNRS, IRIF, F-75013, Paris, France.}
\eaddress{wrigley@irif.fr}

\thanks{
This work has its origin in a question I asked Marco Abbadini at his talk at the Logic Colloquium 2025.  I thank Marco for his stimulating talk as well as the conference organisers for the wonderful event.  I acknowledge the support of Agence Nationale de la Recherche project ANR-23-CE48-0012-01.
}

\title{Existential completions and Herbrand's theorem}


\renewcommand{\phi}{\varphi}
\newcommand{\Doc}[1]{\{{#1}\}\text{-}\mathbf{Doc}}
\newcommand{\DocSites}{\mathbf{DocSites}}
\newcommand{\GeomDoc}{\mathbf{GeomDoc}}
\newcommand{\set}[2]{\{{#1}\,\vert\,{#2}\}}
\newcommand{\cat}{\mathcal{C}}
\newcommand{\dcat}{\mathcal{D}}
\newcommand{\op}{^{\rm op}}
\newcommand{\MLat}{\mathbf{MLat}}
\newcommand{\DLat}{\mathbf{DLat}}
\newcommand{\BA}{\mathbf{BA}}
\newcommand{\Pos}{\mathbf{Pos}}
\newcommand{\Frm}{\mathbf{Frm}}
\newcommand{\Sets}{\mathbf{Set}}
\newcommand{\theory}{\mathbb{T}}
\newcommand{\thdoc}[2]{F_{\{{#2}\}}({#1})}
\newcommand{\Mod}[1]{{#1}\text{-}\mathbf{Mod}}
\newcommand{\Term}{\mathbf{Term}_\Sigma}
\newcommand{\Sh}{\mathrm{Sh}}
\newcommand{\Sub}{\mathrm{Sub}}

\newtheorem{thm}{Theorem}
 
\newtheorem{lem}{Lemma}
\newtheorem{prop}{Proposition}
\newtheorem{coro}{Corollary}

\newtheoremrm{rem}{Remark}
\newtheoremrm{df}{Definition}

\mathrmdef{Hom}
\mathbfdef{Set}

\begin{document}
	
	\maketitle
	\begin{abstract}
		Recently, Abbadini and Guffanti gave an algebraic proof of Herbrand's theorem using a completion for Lawvere doctrines that freely adds existential and universal quantifiers.  A more direct argument can be given by only completing with respect to existential quantifiers.  We construct the free existential completion on a presheaf of distributive lattices, and deduce Herbrand's theorem for coherent logic from the explicit description.  We also discuss the cases involving presheaves of meet-semilattices, due to Trotta, and presheaves of frames.
	\end{abstract}
\section*{Introduction}
\subsection*{Background}
A classical result in proof theory, Herbrand's theorem \cite{herbrand}, asserts that a \emph{universal} first-order theory, i.e.\ one whose axioms avoid any use of existential quantification, can only admit `constructive proofs' of existential statements.  That is to say, if a universal (classical) theory $\theory$ proves that $\top \vdash \exists x.\, \phi(x)$, where $\phi(x)$ is a quantifier-free formula, then there are a finite number of closed terms $t_1, \dots , t_n$ that witness the formula $\phi(x)$, i.e.\ $\theory$ proves $\top \vdash \phi(t_1) \lor \dots \lor \phi(t_n)$ (in \cite[\S 2]{buss2}, this is called the \emph{weak} form of Herbrand's theorem, but it is the only version we will consider in this paper).  Most textbook accounts of Herbrand's theorem employ cut elimination  (see \cite[\S 2.5.1]{buss} for instance), a modern technique that was not available when Herbrand first proved the result.

In a recent work \cite{abba_guff}, Abbadini and Guffanti give an algebraic proof of Herbrand's theorem phrased in the language of \emph{Lawvere doctrines}.  Doctrines {\it\`a la} Lawvere, introduced in \cite{lawvere1,lawvere2}, are a predicate generalisation of the familiar Lindenbaum-Tarski algebras from propositional logic, an alternative to the polyadic algebras of Halmos \cite{halmos} and the {cylindric algebras} of Tarski \cite{tarski}.  A recurring topic within doctrine theory is the study of \emph{completions}: universally adding richer structure to the doctrine.  For example, the \emph{quotient completion} of a doctrine has been extensively studied \cite{quotcomp,quotcompfound,uniexactcomp,maiettirosolinirelating,cioffo}, paralleling the category-theoretic exact completion \cite{carboni,carbonivitale,Carboni_Magno} and the `tripos-to-topos' construction found in \cite{triposoriginal,triposphd,triposretro}.  Recently, more focus has been given to completions that leave the base category of a doctrine unchanged and instead add structure to the fibres: Trotta builds the \emph{existential completion} in \cite{trotta} while the \emph{geometric completion} is considered in \cite{caramello,geom_comp,wrigley_phd}.

Across a pair of papers \cite{abba_guff_2,abba_guff}, Abbadini and Guffanti study the completion of a presheaf of Boolean algebras $P \colon \cat\op \to \BA$ with respect to both universal and existential quantifiers.  They were motivated to develop a doctrinal language in which to phrase problems of quantifier complexity, as sketched in a talk of Gehrke \cite{gehrke}.  The quantifier completion $P^{\forall \exists}$ of $P$, to arbitrary quantifier depth, is proven to exist via an abstract argument in \cite[Theorem 5.2]{abba_guff_2}.  In their follow-up paper \cite{abba_guff}, Abbadini and Guffanti study the completion up to quantifier depth 1, denoted by $P^{\forall\exists}_1$, i.e.\ the free addition of existential and universal formulae to the quantifier-free ones.  To understand the structure of $P^{\forall \exists}_1$, they identify the formulae of quantifier depth 1 as a subdoctrine of $P^{\forall\exists}$ via a semantic condition (their notion of a \emph{universal ultrafilter} 
\cite[Definition 3.3]{abba_guff}).

A careful analysis of the order relation in $P^{\forall\exists}_1$ recovers Herbrand's theorem \cite[Theorem 4.9]{abba_guff}.  The idea is the following: the datum of a universal theory $\theory$ is essentially contained in a presheaf of Boolean algebras -- the doctrine of quantifier-free formulae -- and so its free completion is the doctrine of formulae of quantifier depth 1 modulo $\theory$-provability.  The order relation in the completed doctrine thus describes, algebraically, the statement of Herbrand's theorem.


\subsection*{Our contribution}
The purpose of this paper is to point out that a more direct doctrinal proof of Herbrand's theorem is possible by only completing with respect to existential quantifiers, which in comparison is less involved than the completion with respect to both quantifiers constructed in \cite{abba_guff}.  The first observation to make is that a proof of Herbrand's theorem for \emph{Horn logic} is implicit in the construction of Trotta's existential completion of a primary doctrine \cite{trotta}.  In \cite[Lemma IV.38]{wrigley_phd}, it is shown how Trotta's existential completion can be carved out of the topos-theoretic \emph{geometric completion} found in \cite{geom_comp,caramello,wrigley_phd}, similar to how the exact completion of a regular category is constructed in \cite{lack_exact}.  In a similar fashion (Remark \ref{rem:efface}), the geometric completion instructs us on how to modify Trotta's construction in order to obtain the two novel completions presented in this paper:
\begin{enumerate}
	\item the free existential completion on a doctrine interpreting $\{\land,\lor,=\}$,
	\item and the free existential completion on a doctrine interpreting $\{\land,\bigvee,=\}$.
\end{enumerate}
As a result of the explicit description of these completions we can deduce versions of Herbrand's theorem for coherent logic and geometric logic, with the former subsuming the classical version of Herbrand's theorem (see Remark \ref{rem:morleyisation}).  In contrast, it is not clear how the methods of \cite{abba_guff} generalise beyond the classical setting.
\subsection*{Overview}
We proceed as follows:
\begin{enumerate}[label=(\arabic*)]
	\item In \cref{sec:doctrines}, we recall the theory of Lawvere doctrines and how this acts as an algebraic framework for studying first-order theories, including how existential quantification and the equality predicate are interpreted algebraically.
	\item In \cref{sec:excomp}, we first recall the construction of Trotta's existential completion for presheaves of meet-semilattices, i.e.\ doctrines interpreting only (finite) conjunction, and that this completion respects the interpretation of equality.  Next, we construct the existential completion for presheaves of distributive lattices, and then for presheaves of frames.  For the sake of narrative clarity, and because the proofs closely mirror those in \cite{trotta}, we initially state the universal properties of our existential completions as fact.  We include explicit proofs of the necessary properties in \cref{sec:proofs}.
	\item In \cref{sec:herbrand}, we explain how the explicit descriptions of the existential completions described in the previous section yield versions of Herbrand's theorem for Horn, coherent, classical and geometric logic.
	\item Finally, in \cref{sec:proofs}, we give explicit proofs that our construction does in fact yield the existential completion on a presheaf of distributive lattices.  Explicit proofs that our construction describes the existential completion on a presheaf of frames can be obtained by easily modifying our argument.
\end{enumerate}
\section{Doctrines}\label{sec:doctrines}
Lawvere doctrines, as introduced in \cite{lawvere1,lawvere2}, come in many different flavours depending on the fragment of logic they are intended to interpret.  We will discuss exclusively fragments of \emph{geometric} logic, which include \emph{coherent} logic (also called \emph{positive} logic, for instance, in \cite{benyaacov_poizat}) and therefore subsumes classical logic (\cite[Lemma D1.5.13]{elephant}).  Recall that in geometric logic and its fragments, we use a Gentzen-style sequent calculus, where a sequent $\varphi(\vec{x}) \vdash_{\vec{x}} \psi(\vec{x})$ expresses ``for all $\vec{x}$, if $\phi(\vec{x})$ then $\psi(\vec{x})$''.  For an introduction to the syntax of geometric logic, the reader is directed to \cite[\S D1]{elephant}.
\begin{df}  By a \emph{doctrine} unqualified, we mean a presheaf of posets $P \colon \cat\op \to \Pos$.  Given an arrow $f \colon c' \to c$ in $\cat$, we will write $f^\ast \colon P(c) \to P(c')$ rather than $P(f)$.
\begin{enumerate}
	\item A $\{\land\}$-\emph{doctrine} is a functor $P \colon \cat\op \to \MLat$, where $\cat$ is a category with finite products and $\MLat$ denotes the category of meet-semilattices (i.e.\ posets with a top element and binary joins).  In \cite{quotcomp,trotta}, these are called \emph{primary doctrines}.  The elements $x \in P(c)$ should be thought of as ``predicates in context'', and so in a $\{\land\}$-doctrine we can interpret conjunctions of predicates.
	\item A $\{\land,\lor\}$-\emph{doctrine} is a $\{\land\}$-doctrine that factors through the subcategory of distributive lattices $\DLat \subseteq \MLat$, allowing for disjunctions of predicates.
	\item A $\{\land,\bigvee\}$-\emph{doctrine} is a $\{\land\}$-doctrine that factors through the subcategory of frames $\Frm \subseteq \MLat$. 
\end{enumerate}	
\end{df}
\begin{rem}\label{rem:base_has_products}
	Note that, in some other places (e.g.\ \cite{geom_comp,wrigley_phd}), the term \emph{primary doctrine} assumes all finite limits on the \emph{base} category $\cat$.  This is often a convenient assumption, especially when working with the Grothendieck construction of a doctrine (cf.\ Remark \ref{rem:efface}), but it won't be needed for the constructions in this paper. (It is possible to avoid assuming any limits in the base category, even products, by employing \emph{relative} definitions, see \cite[\S 5]{caramello} or \cite[\S 3]{int_loc_morphisms}).
\end{rem}
So far, the doctrines described above only interpret the logical connectives that were already present in propositional logic.  This is why, for example, a $\{\land,\lor\}$-doctrine is merely a functor that takes values in the category giving algebraic semantics for coherent propositional logic, i.e.\ $\DLat$.  The key insight of Lawvere in \cite{lawvere1} was that first-order logical connectives can be captured by structure \emph{between} the fibres.
\begin{df}[Lawvere, \cite{lawvere1,lawvere2}]\label{df:doctrine} Let $P \colon \cat\op \to \MLat$ be a $\{\land\}$-doctrine.
	\begin{enumerate}
		\item We say that $P$ is a $\{\land,=\}$-\emph{doctrine} if it is equipped with an \emph{equality predicate} $\delta_d \in P(d \times d)$, for each $d \in \cat$, such that the map
		\[\pi_d^\ast(-) \land \pi_{d'}^\ast \delta_d \colon P(d' \times d) \to P(d' \times d \times d)\] defines a left adjoint to $(1_{d'} \times \Delta_d)^\ast \colon P(d' \times d \times d) \to P(d' \times d)$ for each $d' \in \cat$, where $\Delta_d \colon d \to d \times d$ denotes the diagonal map.  In \cite{lawvere2}, this is called \emph{elementary structure} on a doctrine.  Similarly, we will say that $P$ is a $\{\land,\lor,=\}$-\emph{doctrine} (respectively, $\{\land,\bigvee,=\}$-\emph{doctrine}) if $P$ is both a $\{\land,=\}$-doctrine and a $\{\land,\lor\}$-doctrine (resp., $\{\land,\bigvee\}$-doctrine). 
		
		\item We say that $P$ is a $\{\exists,\land\}$-\emph{doctrine} or that $P$ is \emph{existential} if, for each product projection $\pi_d \colon d \times c \to c$ in $\cat$, the map $\pi_d^\ast \colon P(c) \to P(d \times c)$ has a left adjoint $\Sigma_d \colon P(d \times c) \to P(c)$ which moreover satisfies:
		\begin{enumerate}
			\item the \emph{Frobenius} condition that these left adjoints are stable under meets, i.e.\ we have that
			\[
		x \land	\Sigma_d y = \Sigma_d( \pi_d^\ast x \land y)
			\]
			for all $x \in P(c)$ and $y \in P(d\times c)$,
			\item and the \emph{Beck-Chevalley} condition that these left adjoints are stable under pullback, i.e.\ for any arrow $f \colon c' \to c$ in $\cat$, the square
			\[
			\begin{tikzcd}
				P(d \times c) \ar{r}{(1_d \times f)^\ast} \ar{d}[']{\Sigma_d} & P(d \times c) \ar{d}{\Sigma'_d} \\
				P(c) \ar{r}{f^\ast} & P(c').
			\end{tikzcd}
			\]
			As before, we will say that $P$ is a $\{\exists,\land,\lor\}$-doctrine if it is both a $\{\exists,\land\}$-doctrine and a $\{\land,\lor\}$-doctrine, etc.
		\end{enumerate}
	\end{enumerate}
\end{df}
\begin{example}
	Let $\cat$ be a lex category (i.e.\ a category with all finite limits), then taking the subobject functor yields a $\{\land,=\}$-doctrine
	\[
	\Sub_\cat \colon \cat\op \to \MLat
	\]
	where the equality predicate at an object $d \in \cat$ is given by the diagonal subobject $\Delta_d \colon d \rightarrowtail d \times d$.  If $\cat$ is a \emph{regular} category, i.e.\ each arrow $f \colon c' \to c$ in $ \cat$ has an image factorisation $c' \twoheadrightarrow \mathrm{im}(f) \rightarrowtail c$, then $\Sub_\cat$ is a $\{\exists,\land,=\}$-doctrine with the left adjoint $\Sigma_d$ to $\pi_d^\ast$ given by sending a subobject $m \colon e \rightarrowtail d \times c $ to the image $\mathrm{im}(\pi_d \circ m)$ \cite[Lemma A1.3.1]{elephant}.  Similarly, if $\cat$ is a \emph{coherent} category or a \emph{geometric} category (see \cite[\S A1.4]{elephant}) then $\Sub_\cat$ is a $\{\exists,\land,\lor,=\}$-doctrine (respectively, $\{\exists,\land,\bigvee,=\}$-doctrine).  In particular, the \emph{powerset doctrine} $\mathscr{P} \colon \Sets\op \to \Frm$ is a $\{\exists,\land,\bigvee,=\}$-doctrine.
\end{example}
\begin{example}[The syntactic doctrine of a first-order theory]\label{ex:thdoc}
	Let $\theory$ be first-order theory over a signature $\Sigma$.  By $\Term$ we denote the category of \emph{terms} over $\Sigma$, i.e.\ the category where:
	\begin{enumerate}
		\item the objects of $\Term$ are finite contexts $(x_1 , \dots , x_n)$ of variables,
		\item and the morphisms $(x_1, \dots , x_n) \to (y_1, \dots , y_m)$ are tuples of terms with free variables $(t_1(x_1, \dots , x_n) , \dots , t_m(x_1, \dots , x_n))$, where the term $t_i(x_1, \dots , x_n)$ is of the same sort as $y_i$.  The composition of morphisms is given by substitution of terms.
	\end{enumerate}
	Note that the category $\Term$ has finite products: the terminal object is given by the empty context $\emptyset$, while the product is $(x_1,\dots,x_n) \times (y_1 , \dots , y_m) = (x_1, \dots , x_n, y_1, \dots , y_m)$.
	
	We define a $\{\land,=\}$-doctrine $\thdoc{\theory}{\land} \colon \Term \to \MLat$ by taking $\thdoc{\theory}{\land,=}(x_1, \dots ,x_n)$ to be the poset whose elements are $\theory$-provable equivalence classes of $\{\land,=\}$-\emph{formulae}, or \emph{Horn formulae}, in context $(x_1,\dots,x_n)$, i.e.\ the formulae of the form $\chi_1 \land \dots \land \chi_n$ where each $\chi_i$ is either:
	\begin{enumerate}
		\item simply \emph{truth} $\top$,
		\item an instance of equality $s = t$, where $s$ and $t$ are both terms over $\Sigma$ of the same sort with free variables $(x_1, \dots ,x_n)$,
		\item or an instance of a relation $R(t_1, \dots , t_k)$, where $R$ is a relation symbol in $\Sigma $ and each $t_\ell$ is a term (of the appropriate sort) with free variables $(x_1, \dots , x_n)$.
	\end{enumerate}  
	We order these formulae by provability modulo the theory $\theory$, i.e.\ if $\theory$ proves the sequent $\psi \vdash_{\vec{x}} \phi$ then $\psi \leqslant \phi$.  For a morphism $(t_1, \dots , t_m) \colon (x_1, \dots , x_n) \to (y_1 , \dots , y_m)$ in $\Term$, the transition map 
	\[
	\thdoc{\theory}{\land,=}(t_1, \dots , t_m) \colon \thdoc{\theory}{\land,=}(y_1,\dots, y_m) \to \thdoc{\theory}{\land,=}(x_1, \dots , x_n)
	\]
	acts by sending a $\{\land,=\}$-formula $\phi(\vec{y})$ to the result $\phi(\vec{t}/\vec{y})$ of substituting $t_i(\vec{x})$ for $y_i$ (since contexts are assumed to be disjoint, there are no problems with simultaneous substitution).  The equality predicate, at a context $(x_1,\dots,x_n)$, is unsurprisingly given by
	\[
	(x_1 = x'_1) \land \dots \land (x_n = x'_n) \in \thdoc{\theory}{\land,=}(x_1, \dots , x_n, x'_1, \dots , x'_n).
	\]
	If we also include formulae constructed using finitary disjunction $\lor$ or existential quantification $\exists$, we obtain a $\{\land,\lor,=\}$-doctrine that we denote $\thdoc{\theory}{\land,\lor,=}$ (respectively, a $\{\exists,\land,=\}$-doctrine $\thdoc{\theory}{\exists,\land,=}$), and similarly for infinitary disjunction.  More details can be found in \cite{seely}.
	
	If $\theory$ is a propositional coherent theory, i.e.\ the signature $\Sigma$ has no sorts, then $\Term$ is the one-point category $\ast$ and $\thdoc{\theory}{\land,\lor}(\ast)$ is the usual Lindenbaum-Tarski algebra of coherent formulae for $\theory$.
\end{example}
\begin{df}\label{df:docmorph}
	Let $P \colon \cat\op \to \Pos$ and $Q \colon \dcat\op \to \Pos$ be doctrines.  A morphisms of doctrines consists of a pair $(F,\alpha)$ of a functor $F \colon \cat \to \dcat$ and a natural transformation $\alpha \colon P \Rightarrow Q \circ F\op$.  We make the following distinctions.
	\begin{enumerate}
		\item If $P$ and $Q$ are $\{\land\}$-doctrines, then $(F,\alpha)$ is a \emph{morphism of $\{\land\}$-doctrines} if $F$ preserves finite products and $\alpha_c \colon P(c) \to Q(Fc)$ is a homomorphism of meet semi-lattices, for each $c \in \cat$.
		\item If $P$ and $Q$ are $\{\land,\lor\}$-doctrines, then we say that a morphism of $\{\land\}$-doctrines $(F,\alpha) \colon P \to Q$ is a \emph{morphism of $\{\land,\lor\}$-doctrines} (respectively, a \emph{morphism of $\{\land,\bigvee\}$-doctrines}) if $\alpha_c \colon P(c) \to Q(Fc)$ is a homomorphism of distributive lattices (resp., frame homomorphism).
		\item If $P$ and $Q$ are $\{\land,=\}$-doctrines, then we say that a morphism of $\{\land\}$-doctrines $(F,\alpha) \colon P \to Q$ is a \emph{morphism of $\{\land,=\}$-doctrines} if the equality predicate is preserved, i.e.\ $ \alpha_{d \times d}(\delta^P_d) = \delta^Q_{Fd}$ for each $ d \in \cat$.
		\item If $P$ and $Q$ are $\{\exists,\land\}$-doctrines, then we say that a morphism of $\{\land\}$-doctrines $(F,\alpha) \colon P \to Q$ is a \emph{morphism of $\{\exists,\land\}$-doctrines} if the left adjoints to $\pi_d^\ast$ are preserved, i.e.\ for each $c, d \in \cat$ the square
		\[
		\begin{tikzcd}
			P( d \times c) \ar{r}{\alpha_{d \times c}} \ar{d}[']{\Sigma^P_d} & Q(F(d \times c)) \ar{d}{\Sigma_{Fd}^Q} \\
			P(c) \ar{r}{\alpha_c} & Q(Fc)
		\end{tikzcd}
		\]
		commutes.
	\end{enumerate}
	Just as before, we will say that $(F,\alpha)$ is a \emph{morphism of $\{\exists,\land,\lor\}$-doctrines} if $(F,\alpha)$ is both a morphism of $\{\exists,\land\}$-doctrines and of $\{\land,\lor\}$-doctrines, etc.
	
	A \emph{transformation} between a pair of morphisms of doctrines $(F,\alpha), (G,\beta) \colon P \rightrightarrows Q$ consists of a natural transformation $\omega \colon F \Rightarrow G$ such that $\alpha_c(x) \leqslant \omega_c^\ast \beta_c(x)$ for all $c \in \cat$ and $x \in P(c)$.  Thus, by taking a class of doctrines from Definition \ref{df:doctrine}, the corresponding class of morphisms from above, and all transformations between these, we obtain the definition of the following 2-categories of doctrines:
	\begin{align*}
		\Doc{\land}, && \Doc{\exists,\land }, && \Doc{\land,=}, && \Doc{\exists,\land,=} , \\
		\Doc{\land,\lor}, && \Doc{\exists,\land,\lor}, && \Doc{\land,\lor,=}, && \Doc{\exists,\land,\lor,=}, \\
		\Doc{\land,\bigvee}, && \Doc{\exists,\land,\bigvee}, && \Doc{\land,\bigvee,=}, && \Doc{\exists,\land,\bigvee,=}.
	\end{align*}
\end{df}
\begin{example}[Kock \& Reyes, \cite{kockreyes}]\label{ex:morphism_is_model}
	Let $\theory$ be a coherent theory, which for convenience we suppose is single-sorted, and let $\thdoc{\theory}{\exists,\land,\lor,=}$ be the associated syntactic doctrine from Example \ref{ex:thdoc}.  A morphism of $\{\exists,\land,\lor,=\}$-doctrines $(\tilde{M},[-]_M) \colon \thdoc{\theory}{\exists,\land,\lor,=} \to \mathscr{P}$ to the powerset doctrine is precisely the data of a model $M \vDash \theory$ in the usual Tarskian sense.  The product-preserving functor $\tilde{M} \colon \Term \to \Sets$ sends a context $(x_1,\dots,x_n) \in \Term$ to the set of $n$-tuples $|M|^n$, while the lattice homomorphism $[-]_M \colon \thdoc{\theory}{\exists,\land,\lor,=} \to \mathscr{P}(|M|^n)$ sends a formula $\phi(x_1,\dots,x_n)$ to its interpretation $[\phi]_M \subseteq |M|^n$.  Similarly, given a pair of doctrine morphisms $(\tilde{M},[-]_M), (\tilde{N},[-]_N) \colon \thdoc{\theory}{\exists,\land,\lor,=} \to \mathscr{P}$, i.e.\ a pair of models, then a transformation $(M,[-]_M) \Rightarrow (N,[-]_N)$ is precisely the datum of a homomorphism of models.
	
	In the same manner, if we relax the domain of where we take semantics for $\theory$ to be any $\{\exists,\land,\lor,=\}$-doctrine $P \colon \cat\op \to \DLat$, we end up with the essentially tautological observation that a $\theory$-model \emph{internal} to $P$ is a morphism of $\{\exists,\land,\lor,=\}$-doctrines $\thdoc{\theory}{\exists,\land,\lor,=} \to P$.  (If we take $P$ to be the subobject doctrine of a topos, we return the notion of semantics for $\theory$ found in \cite[\S D1.2]{elephant}.)  In this sense, the syntactic doctrine of $\theory$ has the universal property that there is a natural isomorphism of pseudo-functors $\Doc{\exists,\land,\lor,=}(\thdoc{\theory}{\exists,\land,\lor,=},-) \cong \Mod{\theory}(-)$ (analogous to the universal property of a \emph{classifying topos}, cf.\ \cite[\S D3.1]{elephant}).
\end{example}
\section{Existential completions}\label{sec:excomp}
An \emph{existential completion} is a left 2-adjoint to the one of the forgetful 2-functors
\begin{align*}
	\Doc{\exists,\land} & \hookrightarrow \Doc{\land}, & \Doc{\exists,\land,=} & \hookrightarrow \Doc{\land,=}, \\
	\Doc{\exists,\land,\lor} & \hookrightarrow \Doc{\land,\lor}, & \Doc{\exists,\land,\lor,=} & \hookrightarrow \Doc{\land,\lor,=}, \\
	\Doc{\exists,\land,\bigvee} & \hookrightarrow \Doc{\land,\bigvee}, & \Doc{\exists,\land,\bigvee,=} & \hookrightarrow \Doc{\land,\bigvee,=},
\end{align*}
as such a left 2-adjoint would describe how to freely add existential quantifiers to a doctrine of the appropriate type.  As we will see in this section, each of these forgetful 2-functors possesses such an existential completion.
\subsection{The existential completion for $\{\land\}\text{-}$doctrines}
First, recall from \cite{trotta} that there is an existential completion for $\{\land\}$-doctrines and $\{\land,=\}$-doctrines:
\begin{thm}[Trotta \cite{trotta}]\label{thm:exist_comp_for_primary}
	The forgetful 2-functor $\Doc{\exists,\land} \hookrightarrow \Doc{\land}$ has a left 2-adjoint $(-)^\exists \colon \Doc{\land} \to \Doc{\exists,\land}$.  
\end{thm}
Let us recall the explicit construction of the existential completion $(-)^\exists \colon \Doc{\land} \to \Doc{\exists,\land}$ in more detail.  Given a $\{\land\}$-doctrine $P \colon \cat\op \to \MLat$, the completion $P^\exists \colon \cat\op \to \MLat$ has the following description:
\begin{enumerate}
	\item For each object $c \in \cat$, let $P^\exists(c)'$ be the preorder consisting of pairs $(d,x)$, where $d$ is an object of $\cat$ and $x \in P(d \times c)$, with the order given by $(d,x) \leqslant (e,y)$ if there is an arrow $r \colon d \times c \to e \times c$ such that the triangle
	\[
	\begin{tikzcd}
		& e \times c \ar{d}{\pi_e} \\
		d \times c \ar{ru}{r} \ar{r}[']{\pi_d} & c
	\end{tikzcd}
	\]
	commutes and $x \leqslant r^\ast y$.  The fibre $P^\exists(c)$ of the existential completion is the posetal reflection of $P^\exists(c)'$, i.e.\ the poset obtained from $P^\exists(c)'$ by quotienting by the equivalence relation $\sim$ where $(d,x) \sim (e,y)$ if and only if $(d,x) \leqslant (e,y)$ and $(e,y) \leqslant (d,x)$.  As is standard, we will abuse notation and not differentiate between the pair $(a,x)$ and its equivalence class.
	
	\item For each arrow $f \colon c' \to c$, the transition map $P^\exists(f)$, which we will henceforth denote by ${f^\exists}^\ast \colon P^\exists(c) \to P^\exists(c')$, sends (the equivalence class of) a pair $(d,x) \in P^\exists(c)$ to the pair $(d,(1_d \times f)^\ast x) \in P^\exists(c')$, where $1_d \times f$ denotes the universally obtained morphism, i.e.\ the map making the square
	\[
	\begin{tikzcd}
		d \times c' \ar{d}[']{\pi'_d} \ar{r}{1_d \times f} & d \times c \ar{d}{\pi_d} \\
		c' \ar{r}{f} & c
	\end{tikzcd}
	\]
	a pullback.
	
	\item The left adjoint $\Sigma_{d}$ to ${\pi_d^\exists}^\ast \colon P^\exists(c) \to P^\exists(d \times c)$ acts by sending a pair $(e,x) \in P^\exists(d \times c)$, recalling that $x \in P(e \times d \times c)$, to the pair $({e \times d},x) \in P^\exists(c)$.
	
	\item The unit $ P \to P^\exists$ of the existential completion is the morphism of $\{\land\}$-doctrines given by the pair $(1_\cat,\eta)$ consisting of the identity functor on $\cat$ and the natural transformation $\eta \colon P \Rightarrow P^\exists$ that sends $x \in P(c)$ to the pair $(1,x) \in P^\exists(c)$, where $1$ denotes the terminal object in $\cat$.  We will abuse notation and denote both the morphism of doctrines and the natural transformation by $\eta$.
\end{enumerate}
\begin{prop}[Trotta, Proposition 6.1 \cite{trotta}]\label{prop:ex_comp_primary_resps_elem}
	The existential completion respects elementary structure, i.e.\ the 2-functor $(-)^\exists$ in Theorem \ref{thm:exist_comp_for_primary} restricts to a left 2-adjoint for the forgetful 2-functor $\Doc{\exists,\land,=} \hookrightarrow \Doc{\land,=}$.
\end{prop}
\subsection{The existential completion for $\{\land,\lor\}\text{-}$doctrines}\label{subsec:ex_comp_coh}
Trotta's existential completion for $\{\land\}$-doctrines can be modified to yield an existential completion for $\{\land,\lor\}$-doctrines.  For now, we state the construction without proof.  Explicit demonstrations of all the necessary properties can be found in \cref{sec:proofs}.
\begin{thm}\label{thm:exist_comp_for_coherent}
	The forgetful 2-functor $\Doc{\exists,\land,\lor} \hookrightarrow \Doc{\land,\lor}$ has a left 2-adjoint $(-)^\exists \colon \Doc{\land,\lor} \to \Doc{\exists,\land,\lor}$.
\end{thm}
We can explicitly describe the 2-functor $(-)^\exists \colon \Doc{\land,\lor} \to \Doc{\exists,\land,\lor}$.  Let $P \colon \cat\op \to \DLat$ be a $\{\land,\lor\}$-doctrine.
\begin{enumerate}
	\item For each $c \in \cat$, let $P^\exists(c)'$ denote the preorder whose elements are finite sets $\{({d_1},x_1), \dots , ({d_n},x_n)\}$ of pairs consisting of objects $d_1 , \dots , d_n \in \cat$ and $x_i \in P(d_i \times c)$.  The order on $P^\exists(c)'$ is given by 
	\[\{({d_1},x_1), \dots , ({d_n},x_n)\} \leqslant \{({e_1},y_1), \dots , ({e_m},y_m)\}\]
	if, for each $i \leqslant n$, there is finite set of arrows \[r_1 \colon d_i \times c \to e_{j_1} \times c, \dots , r_k \colon d_i \times c \to e_{j_k} \times c,\]
	where $j_1, \dots , j_k \leqslant m$, such that each triangle
	\[
	\begin{tikzcd}
		& e_{j_\ell} \times c \ar{d}{\pi_{e_{j_\ell}}} \\
		d_i \times c \ar{ru}{r_\ell} \ar{r}[']{\pi_{d_i}} & c
	\end{tikzcd}
	\]
	commutes and $x_i \leqslant r_1^\ast y_1 \lor \dots \lor r_k^\ast y_k$.  We take $P^\exists(c)$ to be the posetal reflection of $P^\exists(c)'$.

	\item For each arrow $f \colon c' \to c$ in $\cat$, the transition map ${f^\exists}^\ast \colon P^\exists(c) \to P^\exists(c')$ acts by sending the (equivalence class of) an element $\{({d_1},x_1), \dots , ({d_n},x_n)\} \in P^\exists(c)$ to the finite set of pairs $\{({d_1},(1_{d_1}\times f)^\ast x_1), \dots , ({d_n},(1_{d_n}\times f)^\ast x_n)\} \in P^\exists(c')$.

	\item The left adjoint $\Sigma_{d} \colon P^\exists(d \times c) \to P^\exists(c)$ to ${\pi_d^\exists}^\ast \colon P^\exists(c) \to P^\exists(d \times c)$ acts by sending an element $\{({e_1},x_1), \dots , ({e_n},x_n)\} \in P^\exists(d \times c)$ to the finite set of pairs $\{({e_1 \times d},x_1), \dots , ({e_n \times d},x_n)\} \in P^\exists( c)$.
	
	\item The unit $\eta \colon P \to P^\exists$ is given by the natural transformation that sends $x \in P(c)$ to the singleton pair $\{(1,x)\} \in P^\exists(c)$.  As before, we refer to both the unit and this natural transformation by $\eta$.
\end{enumerate}
Just as in Proposition \ref{prop:ex_comp_primary_resps_elem}, the existential completion on $\{\land,\lor\}$-doctrines respects any elementary structure that is present, in that we have the following result.
\begin{prop}\label{prop:univ_prop_coh_and_equality}
	The 2-functor $(-)^\exists$ in Theorem \ref{thm:exist_comp_for_coherent} restricts to a left 2-adjoint to the forgetful functor $\Doc{\exists,\land,\lor,=} \hookrightarrow \Doc{\land,\lor,=}$.
\end{prop}
\begin{example}[cf.\ Proposition 5.7 \& Corollary 5.8 \cite{abba_guff_2}]\label{ex:universal_doctrine}
	Let $\theory$ be a coherent theory and let $\theory_\forall$ denote the theory of \emph{universal} consequences of $\theory$, by which we mean those sequents $\phi \vdash_{\vec{x}} \psi$ that $\theory$ proves where $\phi, \psi$ are \emph{quantifier-free} coherent formulae, i.e.\ formulae constructed using the logical connectives $\{\land,\lor,=\}$.  Such a sequent is `universal' in the sense that it expresses the sentence ``$\forall x. \, \phi \rightarrow \psi$'', where $\phi$ and $\psi$ are quantifier-free.  We claim that there is an isomorphism
	\[\thdoc{\theory}{\land,\lor,=}^\exists \cong \thdoc{\theory_\forall}{\exists,\land,\lor,=}.\]
	Here $\thdoc{\theory}{\land,\lor,=} \colon \Term \to \DLat$ is the presheaf of quantifier-free formulae, modulo $\theory$-provability. 
	
	The isomorphism can be witnessed by observing that both doctrines satisfy the same universal property.  Recall from Example \ref{ex:morphism_is_model} that a morphism of $\{\exists,\land,\lor,=\}$-doctrines $\thdoc{\theory_\forall}{\exists,\land,\lor,=} \to P$ is the same data as a $\theory_\forall$-model internal to $P$.  By Proposition \ref{prop:univ_prop_coh_and_equality}, a morphism of $\{\exists,\land,\lor,=\}$-doctrines $\thdoc{\theory}{\land,\lor,=}^\exists \to P$ bijectively corresponds to a morphism of $\{\land,\lor,=\}$-doctrines $\thdoc{\theory}{\land,\lor,=} \to P$, i.e.\ an interpretation of each quantifier-free formula in $P$ such that if $\theory$ proves that $\phi \vdash_{\vec{x}} \psi$, with $\phi, \psi$ quantifier free, this inequality is respected in the model.   In other words, a morphism of $\{\land,\lor,=\}$-doctrines $\thdoc{\theory}{\land,\lor,=} \to P$ yields a model of $\theory_\forall$.  Thus, we have that
	\[
	\Doc{\exists,\land,\lor,=}(\thdoc{\theory_\forall}{\exists,\land,\lor,=},-) \cong \Mod{\theory_\forall}(-) \cong \Doc{\exists,\land,\lor,=}(\thdoc{\theory}{\land,\lor,=}^\exists,-)
	\]
	and so we deduce the desired isomorphism.
\end{example}
\begin{rem}\label{rem:efface}
	Lest we `efface our tracks', we describe how one would deduce that the existential completion for a presheaf of distributive lattices is given by the above construction by utilising the theory of \emph{geometric completions} for \emph{doctrinal sites} developed across \cite{caramello,geom_comp,wrigley_phd}.  Since this remark is only adding context, it can safely be skipped by the reader.  For convenience, we will assume that the base category of each doctrine mentioned in this remark has all finite limits, and that the functor underlying each morphism of doctrines preserves these finite limits (cf.\ Remark \ref{rem:base_has_products}).  
	
	Given a doctrine $P \colon \cat\op \to \Pos$, let $\cat \rtimes P$ denote the \emph{Grothendieck construction}, i.e.\ the category whose objects are pairs $(c,x)$, where $c \in \cat$ and $x \in P(c)$, and whose morphisms $f \colon (d,y) \to (c,x)$ are those arrows $f \colon d \to c$ in $\cat$ such that $y \leqslant f^\ast x$.  First observe that, if $P$ is a $\{\land,\lor\}$-doctrine, then $\cat \rtimes P$ has all finite limits (by our above assumptions) and can be endowed with the \emph{Grothendieck topology} $J^P_\lor$ whose covering sieves are generated by the families
	\[
	\{1_c \colon (c,x) \to (c,x\lor y), 1_c \colon (c,y) \to (c,x\lor y) \}
	\] 
	(see \cite[\S III]{SGL} for background on Grothendieck topologies).  Informally, this topology contains all the information of the finite joins in $P$.  More formally, given a second $\{\land,\lor\}$-doctrine $Q \colon \dcat\op \to \DLat$, the datum of a morphism of $\{\land,\lor\}$-doctrines is equivalently described by a morphism of fibrations $\cat \rtimes P \to \dcat \rtimes Q$ that preserves finite limits and sends $J^P_\lor$-covering sieves to $J^Q_\lor$-covering sieves.  Thus, the 2-category $\Doc{\land,\lor}$ can be described as a full and faithful subcategory of the category of \emph{doctrinal sites} \cite[\S III.3.2]{wrigley_phd}, the category $\DocSites$ whose:
	\begin{enumerate}
		\item objects are pairs $(P,J)$ of a $\{\land\}$-doctrine $P \colon \cat\op \to \MLat$ and a Grothendieck topology $J$ on $\cat \rtimes P$,
		\item and whose arrows $(P,J) \to (Q,K)$ are those morphisms of $\{\land\}$-doctrines whose corresponding morphisms of fibrations $\cat \rtimes P \to \dcat \rtimes Q$, in addition to preserving finite limits, also send $J$-covering sieves to $K$-covering sieves.
	\end{enumerate}
	(We are omitting the 2-categorical structure for simplicity.  Note also that in \cite{geom_comp,wrigley_phd} the underlying doctrine of a doctrinal site does not need to be a $\{\land\}$-doctrine, another subtlety that we are eliding.) Explicitly, the category $\Doc{\land,\lor}$ is isomorphic to the full subcategory of $\DocSites$ spanned by objects of the form $(P,J^P_\lor)$, for $P$ a $\{\land,\lor\}$-doctrine.

	As shown in \cite[Theorem IV.25]{wrigley_phd}, the \emph{geometric completion} for a doctrinal site yields a geometric completion for $\{\land,\lor\}$-doctrines.  A \emph{geometric doctrine} is a $\{\land,\bigvee\}$-doctrine $P \colon \cat\op \to \Frm$ such that, for each arrow $f \colon c' \to c$ in $\cat$, the transition map $f^\ast \colon P(c) \to P(c')$ has a left adjoint and moreover these left adjoints must satisfy the Frobenius and Beck-Chevalley conditions.  In particular, a geometric doctrine is a $\{\exists,\land,\bigvee,=\}$-doctrine.  A morphism of geometric doctrines is a morphism of $\{\land,\bigvee\}$-doctrines that preserves these left adjoints (in the same sense as in Definition \ref{df:docmorph}).  Thus, \cite[Theorem IV.25]{wrigley_phd} expresses that there is a left adjoint to the forgetful functor
	\[
	\begin{tikzcd}
		{\Doc{\land,\lor}} && \GeomDoc.
		\arrow[""{name=0, anchor=center, inner sep=0}, "\GC_{\{\land,\lor\}}", shift left=2, from=1-1, to=1-3]
		\arrow[""{name=1, anchor=center, inner sep=0}, shift left=2, hook', from=1-3, to=1-1]
		\arrow["\dashv"{anchor=center, rotate=-90}, draw=none, from=0, to=1]
	\end{tikzcd}
	\]
	Moreover, we know how to compute $\GC_{\{\land,\lor\}}(P)$ for a $\{\land,\lor\}$-doctrine $P$.  It is given by the composite
	\[
	\Omega_{\Sh(\cat \rtimes P,J^P_\lor)} \circ t\op \colon \cat\op \to \Frm
	\]
	of the \emph{subobject classifier} of the topos $\Sh(\cat\rtimes P,J^P_\lor)$ with (the opposite of) the functor $t \colon \cat \to \cat \rtimes P$ that sends an object $c \in \cat$ to $(c,\top) \in \cat \rtimes P$ (see \cite[\S 6]{caramello} and \cite[\S IV.1]{wrigley_phd}).  Alternatively, using the explicit description of the subobject classifier for a Grothendieck topos given in \cite[\S III.7]{SGL}, $\GC_{\{\land,\lor\}}(P)(c)$ is the frame of $J^P_\lor$-\emph{closed sieves} on $(c,\top) \in \cat \rtimes P$ ordered by inclusion.

	The geometric completion freely adds both infinite disjunctions and left adjoints to all transition maps, satisfying the Frobenius and Beck-Chevalley conditions.  In particular, $\GC_{\{\land,\lor\}}$ freely adds existential quantifiers, the left adjoints to product projections.  Thus, we can obtain the existential completion of $P$ by carving out the \emph{subdoctrine} of $\GC_{\{\land,\lor\}}(P)$ generated by $P \subseteq \GC_{\{\land,\lor\}}(P)$, existential quantification and finite disjunctions.  To avoid confusion, denote this subdoctrine by $P_\exists \subseteq \GC_{\{\land,\lor\}}(P)$.  Using the description of $\GC_{\{\land,\lor\}}(P)$ in terms of $J^P_\lor$-closed sieves in $\cat \rtimes P$, this is precisely the subdoctrine ${P_\exists} \subseteq \GC_{\{\land,\lor\}}(P)$ whose elements are the $J^P_\lor$-closed sieves $S \in \GC_{\{\land,\lor\}}(P)(c)$ generated by finite sets of arrows $\set{\pi_{e_i} \colon (e_i \times c ,y_i) \to (c,\top)}{i \leqslant n}$, or explicitly sieves of the form:
	\[
	S = \left\{f  \colon (d,x) \to (c,\top) \, \middle\vert \, \begin{matrix}\text{there are } r_1 \colon d \to e_{j_1} \times c, \dots , r_k \colon d \to e_{j_k} \times c \\ \text{such that } f = \pi_{e_{j_\ell}} \circ r_\ell \text{ and } x \leqslant r_1^\ast y_1 \lor \dots \lor r_k^\ast y_k\end{matrix}\right\}.
	\]
	It is now not too difficult to realise an isomorphism between the $\{\exists,\land,\lor\}$-doctrine $P_\exists$ and the doctrine $P^\exists \colon \cat\op \to \DLat$ described above.
	
	To summarise, starting with a $\{\land,\lor\}$-doctrine $P$, we use a Grothendieck topology $J^P_\lor$ to keep track of the finite disjunctions in $P$, we then use a topos-theoretic construction to complete the pair $(P,J^P_\lor)$ to a \emph{geometric} doctrine $\GC_{\{\land,\lor\}}(P)$, and finally prune $\GC_{\{\land,\lor\}}(P)$ appropriately to arrive at the existential completion $P^\exists$.  This is entirely analogous to the method used in \cite{lack_exact} to construct the exact completion of a regular category.
\end{rem}
\subsection{The existential completion for $\{\land,\bigvee\}$-doctrines}\label{subsec:ex_comp_geom}
Similarly, we can incorporate the infinite disjunctions into our narrative by describing an existential completion for $\{\land,\bigvee\}$-doctrines.  The proof that our stated construction does in fact yield the existential completion on a $\{\land,\bigvee\}$-doctrine is easily modified from the argument presented in \cref{sec:proofs}. 
\begin{thm}\label{thm:ex_comp_for_geom}
	The forgetful 2-functor $\Doc{\exists,\land,\bigvee} \hookrightarrow \Doc{\land,\bigvee}$ has a left 2-adjoint $(-)^\exists \colon \Doc{\land,\bigvee} \to \Doc{\exists,\land,\bigvee}$.
\end{thm}
\begin{enumerate}
	\item For each $c \in \cat$, let $P^\exists(c)'$ denote the preorder whose elements are sets $\set{({d_i},x_i)}{i \in I}$ of pairs where $d_i \in \cat$ and $x_i \in P(d \times c)$ for all $i \in I$.  The order on $P^\exists(c)'$ is given by $$\set{({d_i},x_i)}{i \in I} \leqslant \set{({e_j},y_j)}{j \in J}$$ if, for each $i \in I$, there is a set of arrows $\set{r_k \colon d_i \times c \to e_{j_k} \times c}{k \in K}$ such that each $j_k$ is in $J$, each triangle
	\[
	\begin{tikzcd}
		& e_{j_k} \times c \ar{d}{\pi_{e_{j_k}}} \\
		d_i \times c \ar{ru}{r_k} \ar{r}[']{\pi_{d_i}} & c
	\end{tikzcd}
	\]
	commutes and $x_i \leqslant \bigvee_{k \in K} r_k^\ast y_{j_k}$.  We take $P^\exists(c)$ to be the frame obtained by the posetal reflection of $P^\exists(c)'$.
	
	\item For each arrow $f \colon c' \to c$ in $\cat$, the transition map ${f^\exists}^\ast \colon P^\exists(c) \to P^\exists(c')$ sends $\set{({d_i},x_i)}{i \in I} $ to the set $\set{({d_i},(1_{d_i} \times f)^\ast x_i)}{i \in I} \in P^\exists(c')$.
	
	\item The left adjoint $\Sigma_{d} \colon P^\exists(d \times c) \to P^\exists(c)$ to ${\pi_d^\exists}^\ast \colon P^\exists(c) \to P^\exists(d \times c)$ acts by sending $\set{({e_i},x_i)}{i \in I} \in P^\exists(c)$ to $\set{({e_i \times d},x_i)}{i \in I} \in P^\exists(d \times c)$.
	\item The unit $\eta \colon P \to P^\exists$ sends $x \in P(c)$ to the singleton pair $\{(1,x)\} \in P^\exists(c)$.
\end{enumerate}
Once again, this construction respects any elementary structure that is present:
\begin{prop}
	The 2-functor $(-)^\exists$ in Theorem \ref{thm:ex_comp_for_geom} restricts to a left 2-adjoint to the forgetful functor $\Doc{\exists,\land,\bigvee,=} \hookrightarrow \Doc{\land,\bigvee,=}$.
\end{prop}
\begin{rem}
	As in Remark \ref{rem:efface}, we can derive the description of the existential completion on a $\{\land,\bigvee\}$-doctrine by topos-theoretic means using an analogous argument to the one exposited in Remark \ref{rem:efface} -- for a $\{\land,\bigvee\}$-doctrine $P$, we merely replace the Grothendieck topology $J^P_\lor$ described therein with the topology $J^P_{\bigvee}$ whose covering sieves are generated by the families $\set{\textstyle 1_c \colon (c,x_i) \to (c,\bigvee_{i \in I} x_i)}{i \in I}$.
\end{rem}
\begin{rem}[Lax-idempotency]
	The 2-monads arising from the 2-adjunctions described in Theorem \ref{thm:exist_comp_for_primary}, Theorem \ref{thm:exist_comp_for_coherent} and Theorem \ref{thm:ex_comp_for_geom} are all \emph{lax-idempotent} (also called $KZ$-monads) in the sense of \cite[Definition 1.1]{laxidem}.  The lax-idempotency of the existential completion on $\{\land\}$-doctrines is proven in \cite[Theorem 5.6]{trotta}, while the lax-idempotency for the existential completions on $\{\land,\lor\}$-doctrines and $\{\land,\bigvee\}$-doctrines are consequences of \cite[Corollary IV.28 \& Corollary IV.44]{wrigley_phd}.  Lax-idempotent 2-monads are \emph{property-like} in the sense of \cite{prop_like} in that an object can carry an essentially unique algebra structure, if one exists, e.g.\ a $\{\land\}$-doctrine $P$ can interpret existential quantification in at most one way.
\end{rem}
\section{Herbrand-style theorems}\label{sec:herbrand}
Thus, we arrive at the culmination of our narrative: versions of Herbrand's theorem for Horn, coherent, classical and geometric logic.  As a particular case of Example \ref{ex:universal_doctrine}, given a \emph{universal} theory $\theory$ (as defined below), the free existential completion on the syntactic doctrine of quantifier-free formulae associated with $\theory$ is the syntactic doctrine of \emph{existential} formulae associated with $\theory$.  Herbrand-style theorems can therefore be deduced by simply reading off the definition of the order relation in the appropriate existential completion.
\begin{df}  We distinguish between the following kinds of theories:
	\begin{enumerate}
		\item A theory $\theory$ is a \emph{universal Horn} theory if, for each axiom $\phi \vdash_{\vec{x}} \psi$ of $\theory$, both $\phi$ and $\psi$ are quantifier-free \emph{Horn formulae}, i.e.\ they are constructed using the connectives $\{\land,=\}$.
		\item A theory $\theory$ is a \emph{universal coherent} theory if, for each axiom $\phi \vdash_{\vec{x}} \psi$ of $\theory$, both $\phi$ and $\psi$ are quantifier-free \emph{coherent formulae}, i.e.\ they are constructed using the connectives $\{\land,\lor,=\}$.
		\item A theory $\theory$ is a \emph{universal geometric} theory if, for each axiom $\phi \vdash_{\vec{x}} \psi$ of $\theory$, both $\phi$ and $\psi$ are quantifier-free \emph{geometric formulae}, i.e.\ they are constructed using the connectives $\{\land,\bigvee,=\}$.
	\end{enumerate}
\end{df}
\begin{coro}[Herbrand \cite{herbrand}]\label{coro:herbrand}
	\begin{enumerate}
		\item 	If $\theory$ is a universal Horn theory that proves $\top \vdash \exists x.\, \phi(x)$ where $\phi(x)$ is a quantifier-free Horn formula, then there is a closed term $t$ such that $\theory$ proves $\top \vdash \phi(t)$.
		\item If $\theory$ is a universal coherent theory that proves $\top \vdash \exists x.\, \phi(x)$ where $\phi(x)$ is a quantifier-free coherent formula, then there is a finite set of closed terms $t_1, \dots , t_n$ such that $\theory$ proves $\top \vdash \phi(t_1) \lor \dots \lor \phi(t_n)$.
		\item If $\theory$ is a universal geometric theory that proves $\top \vdash \exists x.\, \phi(x)$ where $\phi(x)$ is a quantifier-free geometric formula, then there is a set of closed terms $\set{ t_i }{i \in I}$ such that $\theory$ proves $\top \vdash \bigvee_{i \in I} \phi(t_i)$.
	\end{enumerate}
\end{coro}
\proof
We spell out the details for a universal coherent theory $\theory$, the other cases being analogous.  As $\theory$ is a universal theory, it is equivalent to its theory of universal consequences.  Thus, by Example \ref{ex:universal_doctrine}, we have that $\thdoc{\theory}{\land,\lor,=}^\exists \cong \thdoc{\theory}{\exists,\land,\lor,=}$.  Under this isomorphism, the quantifier free formula $\phi(x) \in \thdoc{\theory}{\exists,\land,\lor}$ is identified with $\eta_x(\phi(x)) = \{(\emptyset,\phi(x))\} \in \thdoc{\theory}{\land,\lor,=}^\exists$ (recall that $\emptyset$ describes the terminal object of $\Term$).  Therefore, $\theory$ proves the sequent $\top \vdash \exists x. \, \phi(x)$  if and only if there is an inequality
\[
\eta_\emptyset(\top) = \{(\emptyset,\top)\} \leqslant \Sigma_x \{(\emptyset,\phi(x))\} = \{(x,\phi(x))\}
\]
(note that $\{(\emptyset,\top)\}$ describes the top element of $\thdoc{\theory}{\land,\lor,=}^\exists(\emptyset)$, see Lemma \ref{lem:pex_is_dlat}).  Recall from \cref{subsec:ex_comp_coh} that this inequality holds if and only if there is a finite family of arrows $\set{t_i \colon \emptyset \to x}{i \leqslant n}$, i.e.\ a finite set of closed terms, such that
\[
\top \vdash t_1^\ast \phi(x) \lor \dots \lor t_n^\ast \phi(x) = \phi(t_1) \lor \dots \lor \phi(t_n).
\] 
\endproof
\begin{rem}
	Of course, Corollary \ref{coro:herbrand} only describes the simplest instance of Herbrand's theorem where we consider a sequent of the form $\top \vdash \exists x.\, \phi(x)$.  As it is sufficiently clear how to calculate the statement from the description of the existential completion, we leave it as an exercise to the reader to write out explicitly what it means for a coherent theory to prove a more complicated sequent such as $\phi(\vec{x}) \land \exists \vec{y}.\, \psi(\vec{x},\vec{y}) \vdash_{\vec{x}} \chi(\vec{x}) \land \exists \vec{z}. \, \xi(\vec{x},\vec{z})$, for $\phi,\psi, \chi, \xi$ quantifier-free coherent formulae.
\end{rem}
\begin{rem}[Herbrand's theorem for classical logic]\label{rem:morleyisation}
Not much must be changed if we wish to incorporate classical logic into our doctrinal account of Herbrand's theorem.  First note that in classical first-order logic the sequent $\phi \vdash_{\vec{x}} \psi$ is equivalent to the sequent $\top \vdash_{\vec{x}} \neg \phi \lor \psi$, and so we can use a Hilbert-style calculus if we prefer.  By a \emph{universal classical} theory $\theory$ we mean a theory of classical first-order logic whose axioms are all of the form $\top \vdash_{\vec{x}} \phi$ where $\phi$ is a quantifier-free classical formula, i.e.\ constructed using the connectives $\{\land,\lor,=,\neg\}$.  Let $\thdoc{\theory}{\land,\lor,=,\neg} \colon \Term \to \BA \subseteq \DLat$ denote the doctrine of quantifier-free classical formulae modulo $\theory$-provability, taking values in the subcategory $\BA$ of Boolean algebras and their homomorphisms.  By an argument identical to that in Example \ref{ex:universal_doctrine}, the existential completion $\thdoc{\theory}{\land,\lor,=,\neg}^\exists$ is the doctrine of \emph{existential formulae}, i.e.\ formulae of the form $\exists \vec{x}. \, \phi(\vec{x},\vec{y})$ where $\phi(\vec{x},\vec{y})$ is a quantifier-free classical formula, ordered by $\theory$-provability.  Thus, just as in Corollary \ref{coro:herbrand}, a universal classical theory $\theory$ proves $\top \vdash \exists x. \, \phi(x)$, where $\phi(x)$ is a quantifier-free classical formula, if and only if there is a finite set of closed terms $t_1, \dots , t_n$ such that $\theory$ proves $\top \vdash \phi(t_1) \lor \dots \lor \phi(t_n)$.

Alternatively, we could substitute the universal classical theory $\theory$ with an equivalent universal coherent theory $\theory'$ by replacing each instance of a negated (quantifier-free) coherent formula $\phi(\vec{x})$ in the axiomatisation of $\theory$ with a new relation symbol $(\neg \phi)(\vec{x})$ satisfying the sequents $\phi(\vec{x}) \land (\neg \phi)(\vec{x}) \vdash_{\vec{x}} \bot$ and $\top \vdash_{\vec{x}} \phi(\vec{x}) \lor (\neg \phi)(\vec{x}) $ (cf.\ the \emph{Morleyisation} of a classical first-order theory \cite[Lemma D1.5.13]{elephant}), to which we can directly apply the statement of Corollary \ref{coro:herbrand}.

By manipulating our sequents using classical first-order logic, we are not restricted to statements involving only existential quantification.  For example, we might ask (as the opening of \cite{abba_guff} does) when a universal classical theory $\theory$ proves $\forall x. \, \phi(x) \vdash \forall y. \, \psi(y)$ where $\phi$ and $\psi$ are quantifier-free classical formulae.  Using classical logic, we can mutate the sequent $\forall x. \, \phi(x) \vdash \forall y. \, \psi(y)$ into a form which only mentions existential quantification and to which we can apply the existential completion of \cref{subsec:ex_comp_coh}, namely we have that the following sequents are provably equivalent:
\begin{align*}
 \forall x. \, \phi(x) & \vdash \forall y. \, \psi(y) \\
		\iff \top & \vdash \neg \forall x. \, \phi(x) \lor \forall y. \, \psi(y),  \\
		\iff \top & \vdash \exists x. \, \neg \phi(x) \lor \forall y. \, \psi(y), \\
		\iff \top  & \vdash \forall y. \, ( \exists x. \, \neg \phi(x) \lor \psi(y)), \\
		\iff \top & \vdash_y \exists x. \, \neg \phi(x) \lor \psi(y),
\end{align*}
where we have used that $\forall x. \, \phi(x)  \equiv \neg \exists x. \,  \neg \phi(x) $ and the `dual' of the Frobenius condition for universal quantification, which is valid in classical logic (see \cite[Remark D1.3.9]{elephant}).  Therefore, just as in Corollary \ref{coro:herbrand}, the universal theory $\theory$ proves the sequent $\top \vdash_y \exists x. \, \neg \phi(x) \lor \psi(y)$ if and only if there is an inequality 
\begin{align*}
	\eta_y(\top) = \{(\emptyset,\top)\} & \leqslant \Sigma_x \pi_y^\ast \eta_x(\neg \phi(x)) \lor \eta_y(\psi(y)) = \Sigma_x \eta_{x,y}(\neg \phi(x)) \lor \eta_y(\psi(y)) \\
	& = \{(x,\neg\phi(x)),(\emptyset,\psi(y))\}
\end{align*}
in the poset $\thdoc{\theory}{\land,\lor,=,\neg}^\exists(y)$, which we recall occurs if and only if there is a finite set of arrows $t_1 \colon y \to x \times y, \dots , t_n \colon y \to x \times y$ in $\Term$ satisfying the commutativity condition $\pi_x \circ t_i = 1_y$, in other words a finite set of terms $t_1(y) , \dots , t_n(y)$ of type $x$ with free variable $y$, for which $\theory$ proves the sequent $\top \vdash_y \neg \phi(t_1(y)) \lor \dots \lor \neg \phi(t_n(y)) \lor \psi(y)$.  Once again, we can use classical logic to manipulate the former sequent and deduce that the following are provably equivalent:
\begin{align*}
\top	& \vdash_y \neg \phi(t_1(y)) \lor \dots \lor \neg \phi(t_n(y)) \lor \psi(y) \\
	\iff \neg ( \neg \phi(t_1(y)) \lor \dots \lor \neg \phi(t_n(y))) & \vdash_y \psi(y), \\
	\iff \quad\quad\ \ \,  \phi(t_1(y)) \land \dots \land \phi(t_n(y)) & \vdash_y \psi(y).
\end{align*}
Hence, we conclude that $\theory$ proves $\forall x. \, \phi(x) \vdash \forall y. \, \psi(y)$ if and only if there are a finite set of terms $t_1(y), \dots , t_n(y)$ such that $\theory$ proves $\phi(t_1(y)) \land \dots \land \phi(t_n(y)) \vdash_y \psi(y)$ (cf.\ Theorem 4.7 \cite{abba_guff}).
\end{rem}
\section{The proofs for the existential completion}\label{sec:proofs}
We now return to giving a full proof that the construction described in \cref{subsec:ex_comp_coh} does indeed yield the existential completion of a $\{\land,\lor\}$-doctrine.  We could argue via the universal property of the geometric completion and the description of $P^\exists$ as a subdoctrine of the geometric completion (Remark \ref{rem:efface}), but for clarity we instead give a direct argument, similar to that found in \cite{trotta}.  Our account is easily modified to the case of $\{\land,\bigvee\}$-doctrines described in \cref{subsec:ex_comp_geom}.
\subsection{Constructing the existential completion on a $\{\land,\lor\}$-doctrine}
First, for a given $\{\land,\lor\}$-doctrine $P \colon \cat\op \to \DLat$, we show that $P^\exists \colon\cat\op \to \DLat$ is a well-defined $\{\exists,\land,\lor\}$-doctrine.
\begin{lem}\label{lem:pex_is_dlat}
	For each $c \in \cat$, the poset $P^\exists(c)$ is a distributive lattice.
\end{lem}
\proof
It is easily checked that the top and bottom elements of $P^\exists(c)$ are given by $\{(1,\top_c)\}$ and $\emptyset$ respectively, and that the join of the pair $\{({d_1},x_1), \dots , ({d_n},x_n)\}$ and $\{({e_1},y_1), \dots , ({e_m},y_m)\}$ is given by the union
\[
\{({d_1},x_1), \dots , ({d_n},x_n),({e_1},y_1), \dots , ({e_m},y_m)\}.
\]
The meet of $\{({d_1},x_1), \dots , ({d_n},x_n)\}$ and $\{({e_1},y_1), \dots , ({e_m},y_m)\}$ is more involved: it is given by $\set{({d_i \times e_j},{\pi'}_{e_j}^\ast x_i \land {\pi'}_{d_i}^\ast y_j)}{i \leqslant n, j \leqslant m}$, where ${\pi'}_{e_j}$ and ${\pi'}_{d_i}$ are the product projections
\[
\begin{tikzcd}
	d_i \times e_j \times c \ar{r}{{\pi'}_{d_i}} \ar{d}[']{{\pi'}_{e_j}} & e_j \times c \ar{d}{\pi_{e_j}} \\
	d_i \times c \ar{r}{\pi_{d_i}} & c.
\end{tikzcd}
\]
Note that this square is a pullback.  To see why this describes the meet, first note that the families of arrows $\set{{\pi'}_{e_j} \colon d_i \times e_j \times c \to d_i \times c }{j \leqslant m}$ and $\set{{\pi'}_{d_i} \colon d_i \times e_j \times c \to e_j \times c}{i \leqslant n}$ witness the inequalities
\[
\set{({d_i \times e_j},{\pi'}_{e_j}^\ast x_i \land {\pi'}_{d_i}^\ast y_j)}{i \leqslant n, j \leqslant m} \leqslant \set{({d_i},x_i)}{i \leqslant n}, \set{({e_j},y_j)}{j \leqslant m}.
\]
Now suppose that
\[
\{({b_1},z_1), \dots , ({b_p},z_p)\} \leqslant \set{({d_i},x_i)}{i \leqslant n}, \set{({e_j},y_j)}{j \leqslant m},
\]
and so, for each $q \leqslant p$, there are finite sets of arrows $r_1 \colon b_q \times c \to d_{i_1} \times c, \dots,  r_k \colon b_q \times c \to d_{i_k} \times c$ and $r'_1 \colon b_q \times c \to e_{j_1} \times c, \dots , r'_{k'} \colon b_q \times c \to e_{j_{k'}} \times c$ such that,  for every $\ell \leqslant k, \ell' \leqslant k'$, the triangles
\[
\begin{tikzcd}
	& d_{i_\ell} \times c \ar{d}{\pi_{d_{i_\ell}}} \\
	b_q \times c \ar{r}[']{\pi_{b_q}} \ar{ru}{r_\ell} & c, 
\end{tikzcd}
\quad
\begin{tikzcd}
	& e_{j_{\ell'}} \times c \ar{d}{\pi_{e_{j_{\ell'}}}} \\
	b_q \times c \ar{r}[']{\pi_{b_q}} \ar{ru}{r'_{\ell'}} & c, 
\end{tikzcd}
\]
commute and $z_q \leqslant r_1^\ast x_{i_1} \lor \dots \lor r_k^\ast x_{i_k}$ and $z_q \leqslant {r'}_1^\ast y_{j_1} \lor \dots \lor {r'}_{k'}^\ast y_{j_{k'}}$.  Thus, we have that
\begin{align*}
	z_q \leqslant & ( r_1^\ast x_{i_1} \lor \dots \lor r_k^\ast x_{i_k} ) \land ( {r'}_1^\ast y_{j_1} \lor \dots \lor {r'}_{k'}^\ast y_{j_{k'}}), \\
	= & \textstyle \bigvee_{\ell \leqslant k, \ell' \leqslant k'} r_\ell^\ast x_{i_\ell} \land {r'}^\ast_{\ell'} y_{j_{\ell'}}, \\ 
	= & \textstyle \bigvee_{\ell \leqslant k, \ell' \leqslant k'} (r_\ell \times r'_{\ell'})^\ast {\pi'}_{e_{j_{\ell'}}}^\ast x_{i_\ell} \land (r_\ell \times r'_{\ell'})^\ast {\pi'}_{d_\ell}^\ast  y_{j_{\ell'}}, \\
	= & \textstyle \bigvee_{\ell \leqslant k, \ell' \leqslant k'} (r_\ell \times r'_{\ell'})^\ast ({\pi'}_{e_{j_{\ell'}}}^\ast x_{i_\ell} \land  {\pi'}_{d_\ell}^\ast  y_{j_{\ell'}}),
\end{align*}
where $r_\ell \times r'_{\ell'}$ denotes the universally induced arrow
\[
\begin{tikzcd}
	b_q \times c \ar[dashed]{rd}{r_\ell \times r'_{\ell'}} \ar[bend left = 2em]{rrd}{r'_{\ell'}} \ar[bend right]{rdd}[']{r_\ell} \\
	& d_{i_\ell} \times e_{j_{\ell'}} \times c \ar{r}{{\pi'}_{d_{i_\ell}}} \ar{d}[']{{\pi'}_{e_{j_{\ell'}}}} & e_{j_{\ell'}} \times c \ar{d}{\pi_{e_{j_{\ell'}}}} \\
	& d_{i_\ell} \times c \ar{r}{\pi_{d_{i_\ell}}} & c.
\end{tikzcd}
\]
Thus, we have that $\{({b_1},z_1), \dots , ({b_p},z_p)\} \leqslant \set{({d_i \times e_j},{\pi'}_{e_j}^\ast x_i \land {\pi'}_{d_i}^\ast y_j)}{i \leqslant n, j \leqslant m}$ as desired.  Finally, joins and meets distribute via the direct calculation:
\begin{align*}
	& ( \set{({d_i},x_i)}{i \leqslant n} \lor \set{({e_j},y_j)}{j \leqslant m}) \land \set{({b_q},z_q)}{q \leqslant p}, \\
	= \ & \set{({d_i},x_i), ({e_j},y_j)}{i \leqslant n, j \leqslant m} \land \set{({b_q},z_q)}{q \leqslant p}, \\
	= \ & \set{({d_i \times b_q},{\pi'}^\ast_{b_q}x_i \land {\pi'}^\ast_{d_i}z_q),({e_j \times b_q},{\pi'}^\ast_{b_q}y_j \land {\pi'}^\ast_{e_j}z_q)}{i \leqslant n, j \leqslant m, q \leqslant p}, \\
	= \ & \set{({d_i \times b_q},{\pi'}^\ast_{b_q}x_i \land {\pi'}^\ast_{d_i}z_q)}{i \leqslant n, q \leqslant p} \lor \set{({e_j \times b_q},{\pi'}^\ast_{b_q}y_j \land {\pi'}^\ast_{e_j}z_q)}{ j \leqslant m, q \leqslant p}, \\
	= \ & (\set{({d_i},x_i)}{i \leqslant n}  \land \set{({b_q},z_q)}{q \leqslant p}) \lor ( \set{({e_j},y_j)}{j \leqslant m} \land \set{({b_q},z_q)}{q \leqslant p})
\end{align*}
\endproof
\begin{lem}
	The function ${f^\exists}^\ast$ is a lattice homomorphism for each arrow in $\cat$.
\end{lem}
\proof
Clearly ${f^\exists}^\ast$ preserves joins since these are just unions.  Meets are preserved since we have that:
\begin{align*}
	& {f^\exists}^\ast ( \set{({d_i},x_i)}{i \leqslant n} \land \set{({e_j},y_j)}{j \leqslant m}) \\
	 = \ & {f^\exists}^\ast \set{({d_i \times e_j}, \pi_{e_j}^\ast x_i \land \pi_{d_i}^\ast y_j)}{i \leqslant n, j \leqslant m} \\
	= \ & \set{({d_i \times e_j}, (1_{d_i \times e_j} \times f)^\ast( {\pi}_{e_j}^\ast x_i \land  {\pi}_{d_i}^\ast y_j))}{i \leqslant n, j \leqslant m}, \\
	= \ & \set{({d_i \times e_j}, (1_{d_i \times e_j} \times f)^\ast {\pi}_{e_j}^\ast x_i \land (1_{d_i \times e_j} \times f)^\ast {\pi}_{d_i}^\ast y_j)}{i \leqslant n, j \leqslant m}, \\
	= \ & \set{({d_i \times e_j}, {\pi'}_{e_j}^\ast (1_{d_i} \times f)^\ast x_i \land {\pi'}_{d_i}^\ast (1_{e_j} \times f)^\ast y_j)}{i \leqslant n, j \leqslant m}, \\
	= \ & {f^\exists}^\ast \set{({d_i},x_i)}{i \leqslant n} \land {f^\exists}^\ast \set{({e_j},y_j)}{j \leqslant m},
\end{align*} 
where we have used the commutativity of the two squares
\[
\begin{tikzcd}
	d_i \times e_j \times c' \ar{r}{1_{d_i \times e_j} \times f} \ar{d}[']{\pi'_{e_j}} & d_i \times e_j \times c \ar{d}{\pi_{e_j}} \\
	d_i \times c' \ar{r}{1_{d_i} \times f} & d_i \times c,
\end{tikzcd} \quad
\begin{tikzcd}
	d_i \times e_j \times c' \ar{r}{1_{d_i \times e_j} \times f} \ar{d}[']{\pi'_{d_i}} & d_i \times e_j \times c \ar{d}{\pi_{d_i}} \\
	e_j \times c' \ar{r}{1_{e_j} \times f} & e_j \times c.
\end{tikzcd}
\]
\endproof
\begin{lem}
	There is an adjunction $\Sigma_d \dashv {\pi_d^\exists}^\ast$ for each $d \in \cat$.
\end{lem}
\proof
Let $\{({e_1},x_1), \dots , ({e_n},x_n)\}$ and $\{({b_1},y_1), \dots , ({b_q},y_q) \}$ be elements of $P^\exists(d \times c)$ and $P^\exists(c)$ respectively.  We have that
\[
\Sigma_d (\{({e_1},x_1), \dots , ({e_n},x_n)\}) 
=  \{({e_1 \times d},x_1), \dots , ({e_n \times d},x_n)\} 
\leqslant  \{({b_1},y_1), \dots , ({b_q},y_q) \}
\]
if and only if, for each $i \leqslant n$, there is a finite set of arrows
\[r_1 \colon e_i \times d \times c \to b_{j_1} \times c, \dots , r_k \colon e_i \times d \times c \to b_{j_k} \times c,\]
where $j_1 , \dots , j_k \leqslant q$, such that each triangle
\[
\begin{tikzcd}
	& b_{j_\ell} \times c \ar{d}{\pi_{b_{j_\ell}}} \\
	e_i \times d \times c \ar{ru}{r_\ell} \ar{r}[']{\pi_{e_i \times d}} & c
\end{tikzcd}
\]
commutes and $x_i \leqslant r_1^\ast y_{j_1} \lor \dots \lor r_k^\ast y_{j_k}$.  Meanwhile, we have that
\begin{align*}
	\{({e_1},x_1), \dots , ({e_n},x_n)\} & \leqslant  {\pi_d^\exists}^\ast ( \{({b_1},y_1), \dots , ({b_q},y_q) \}) \\
	& =  \{({b_1},(1_{b_1} \times \pi_d)^\ast y_1, \dots , ({b_q},(1_{b_q} \times \pi_d)^\ast y_q))\}
\end{align*}
if and only if, for each $i \leqslant n$, there is a finite set of arrows 
\[r'_1 \colon e_i \times d \times c \to b_{j'_1} \times d \times c, \dots, r'_{k'} \colon e_i \times d \times c \to b_{j'_{k'}} \times d \times c,\]
where $j'_1 , \dots , j'_k \leqslant q$, such that each triangle
\[
\begin{tikzcd}
	& b_{j'_\ell} \times d \times c \ar{d}{\pi'_{b_{j'_\ell}}} \\
	e_i \times d \times c \ar{ru}{r'_\ell} \ar{r}[']{\pi'_{e_i \times d}} & c
\end{tikzcd}
\]
commutes and $x_i \leqslant {r'_1}^\ast (1_{b_{j'_1}} \times \pi_d)^\ast y_{j'_1} \lor \dots \lor {r'_{k'}}^\ast (1_{b_{j'_{k'}}} \times \pi_d)^\ast y_{j'_{k'}}$.

Thus, given a family $ r_1 \colon e_i \times d \times c \to b_{j_1} \times c, \dots , r_k \colon e_i \times d \times c \to b_{j_k} \times c$, for each $i \leqslant n$, witnessing that $ \Sigma_d (\{({e_1},x_1), \dots , ({e_n},x_n)\}) \leqslant  \{({b_1},y_1), \dots , ({b_q},y_q) \}$, we obtain a family of arrows witnessing that $\{({e_1},x_1), \dots , ({e_n},x_n)\}  \leqslant  {\pi_d^\exists}^\ast ( \{({b_1},y_1), \dots , ({b_q},y_q) \}) $ by taking, for each $i \leqslant n$, the family of universally induced morphisms
\[
\begin{tikzcd}
	& d \\
	e_i \times d \times c \ar[bend left = 2em]{ru}{\pi_{e_i \times c}} \ar[bend right = 2em]{rd}[']{r_\ell} \ar[dashed]{r}{r'_{j_\ell}} & b_{j_\ell} \times d \times c \ar{u}[']{\pi_{b_{j_\ell} \times c}} \ar{d}{(1_{b_{j_\ell}} \times \pi_d)} \\
	& b_{j_\ell} \times c
\end{tikzcd}
\]
Conversely, given a family $r'_1 \colon e_i \times d \times c \to b_{j'_1} \times d \times c, \dots, r'_{k'} \colon e_i \times d \times c \to b_{j'_{k'}} \times d \times c$, for each $i \leqslant n$, witnessing that $\{({e_1},x_1), \dots , ({e_n},x_n)\}  \leqslant  {\pi_d^\exists}^\ast ( \{(b_1,y_1), \dots , ({b_q},y_q) \}) $, the family of composites 
\[
(1_{b_{j'_1}} \times \pi_d) \circ r'_1 \colon e_i \times e_i \times d \times c \to b_{j'_1} \times c, \dots , (1_{b_{j'_{k'}}} \times \pi_d) \circ r'_{k'} \colon e_i \times d \times c \to b_{j'_{k'} } \times c
\]
witness the inequality $ \Sigma_d (\{({e_1},x_1), \dots , ({e_n},x_n)\}) \leqslant  \{({b_1},y_1), \dots , ({b_q},y_q) \}$.  Therefore, we have that
\begin{align*}
	\Sigma_d (\{({e_1},x_1), \dots , ({e_n},x_n)\}) & \leqslant  \{(b_1,y_1), \dots , (b_q,y_q) \} \\
	\iff  \{({e_1},x_1), \dots , ({e_n},x_n)\}  & \leqslant  {\pi_d^\exists}^\ast ( \{({b_1},y_1), \dots , ({b_q},y_q) \}),
\end{align*}
i.e.\ there is an adjunction $\Sigma_d \dashv {\pi_d^\exists}^\ast$ as desired.
\endproof
\begin{lem}
	The left adjoints $\Sigma_d$, for $d \in \cat$, satisfy the Frobenius and Beck-Chevalley conditions.
\end{lem}
\proof
Given $d \in \cat$, $\set{({e_i},x_i)}{i \leqslant n} \in P^\exists(d \times c) $ and $\set{({b_p},y_p)}{p \leqslant q} \in P^\exists(c)$, for the Frobenius condition we need to show that
\[
\Sigma_d(\set{({e_i},x_i)}{i \leqslant n}) \land \set{({b_p},y_p)}{p \leqslant q} = \Sigma_d( \set{({e_i},x_i)}{i \leqslant n} \land  {\pi_d^\exists}^\ast \set{({b_p},y_p)}{p \leqslant q} ).
\]
This follows by the calculation:
\begin{align*}
	& \Sigma_d(\set{({e_i},x_i)}{i \leqslant n}) \land \set{({b_p},y_p)}{p \leqslant q} \\
	= \ & \set{({e_i \times d} , x_i)}{i \leqslant n} \land \set{({b_p},y_p)}{p \leqslant q}, \\
	= \ & \set{({e_i \times d \times b_p},{\pi}_{b_p}^\ast x_i \land {\pi}_{e_i \times d}^\ast y_p)}{i \leqslant n , p \leqslant q}, \\
	= \ & \Sigma_d (\set{({e_i  \times b_p},{\pi}_{b_p}^\ast x_i \land {\pi}_{e_i \times d}^\ast y_p)}{i \leqslant n , p \leqslant q}), \\
	= \ & \Sigma_d ( \set{({e_i  \times b_p},{\pi}_{b_p}^\ast x_i \land {\pi}_{e_i}^\ast (1_{b_q} \times \pi_d)^\ast y_p)}{i \leqslant n , p \leqslant q} ), \\
	= \ & \Sigma_d(\set{({e_i},x_i)}{i \leqslant n} \land \set{({b_p},(1_{b_q} \times \pi_d)^\ast y_p)}{p \leqslant q}), \\
	= \ & \Sigma_d( \set{({e_i},x_i)}{i \leqslant n} \land  {\pi_d^\exists}^\ast \set{({b_p},y_p)}{p \leqslant q} ).
\end{align*}

For the Beck-Chevalley condition, we need to show that for any arrow $f \colon c' \to c$ in $\cat$, the square
\[
\begin{tikzcd}
	P^\exists(d \times c) \ar{r}{{(1_d \times f)^\exists}^\ast} \ar{d}[']{\Sigma_d} & P^\exists(d \times c') \ar{d}{\Sigma'_d} \\
	P^\exists(c) \ar{r}{{f^\exists}^\ast} & P^\exists(c')
\end{tikzcd}
\]
commutes.  This also holds since, given $\set{({e_i},x_i)}{i \leqslant n} \in P^\exists(d \times c)$, we have that:
\begin{align*}
	\Sigma'_d {(1_d \times f)^\exists}^\ast \set{({e_i},x_i)}{i \leqslant n} = \ & \Sigma'_d \set{({e_i},(1_{e_i} \times 1_d \times f)^\ast x_i)}{i \leqslant n}, \\
	= \  & \set{({e_i \times d},(1_{e_i \times d} \times f)^\ast x_i}{i \leqslant n}, \\
	= \  & {f^\exists}^\ast \set{({e_i \times d},x_i)}{i \leqslant n} =  {f^\exists}^\ast \Sigma_d \set{({e_i},x_i)}{i \leqslant n}.
\end{align*}
\endproof
\begin{coro}
	The doctrine $P^\exists \colon \cat\op \to \DLat$ is a $\{\exists,\land,\lor\}$-doctrine.
\end{coro}
\subsection{The universal property of the existential completion}
We now show that $P^\exists$ does indeed satisfy the universal property of the existential completion as described in Theorem \ref{thm:exist_comp_for_coherent}.
\begin{lem}
	The unit $\eta \colon P \to P^\exists$ is a morphism of $\{\land,\lor\}$-doctrines.
\end{lem}
\proof
We just need to check that $\eta_c \colon P(c) \to P^\exists(c)$, given by $x \mapsto \{(1,x)\}$, preserves finite meets and joins, for each $c \in \cat$.  The preservation of meets is immediate by the construction of meets in $P^\exists(c)$ found in Lemma \ref{lem:pex_is_dlat}.  For joins, we need to show that $\{(1,x\lor y)\} \leqslant \{(1,x),(1,y)\}$ and $\{(1,x),(1,y)\} \leqslant \{(1,x\lor y)\}$.  In both cases, this is immediate with the witnessing morphism(s) being the identity on $c \cong 1 \times c$.
\endproof
\begin{lem}\label{lem:eta_is_inj}
	Moreover, for each $c \in \cat$, the homomorphism $\eta_c \colon P(c) \to P^\exists(c)$ is injective.
\end{lem}
\proof
If $\eta_c(x) = \{(1,x)\} = \{(1,y)\} = \eta_c(y)$, then there are two families of morphisms $\set{r_i \colon 1 \times c \to 1 \times c}{i \leqslant n}$ and $\set{r'_j \colon 1 \times c \to 1 \times c}{j \leqslant m}$ for which $x \leqslant \bigvee_{i \leqslant n} r_i^\ast y$ and $y \leqslant\bigvee_{j \leqslant m} {r'_j}^\ast x$ and making the triangles
\[
\begin{tikzcd}
	& 1 \times c \ar{d} \\
	1 \times c \ar{ru}{r_i} \ar{r} & c,
\end{tikzcd}
\begin{tikzcd}
	& 1 \times c \ar{d} \\
	1 \times c \ar{ru}{r'_j} \ar{r} & c
\end{tikzcd}
\]
both commute.  But this commutativity condition forces each $r_i,r'_j$ to be the identity $1_c$.  Thus, $x \leqslant y $ and $y \leqslant x$ and so $x = y$.
\endproof
\begin{prop}\label{prop:exist_comp_1cat}
	Let $Q \colon \dcat \to \DLat$ be a $\{\exists,\land,\lor\}$-doctrine.  For every morphism $(F,\alpha) \colon P \to Q$ of $\{\land,\lor\}$-doctrines, there is a unique morphism $(F,\alpha^\exists) \colon P^\exists \to Q$ of $\{\exists,\land,\lor\}$-doctrines making the diagram commute
	\[
	\begin{tikzcd}
		P \ar{r}{\eta} \ar{rd}[']{(F,\alpha)} & P^\exists \ar[dashed]{d}{(F,\alpha^\exists)} \\
		& Q.
	\end{tikzcd}
	\]
\end{prop}
\proof
That $(F,\alpha^\exists)$ must be a commuting morphism of $\{\exists,\land,\lor\}$-doctrines leaves only one possible definition:
\[
\alpha_c^\exists(\{({d_1},x_1) , \dots , (d_n,x_n)\}) = \exists^Q_{Fd_1} \alpha_{d_1 \times c}(x_1) \lor \dots \lor \exists^Q_{Fd_n} \alpha_{d_n \times c}(x_n),
\]
where $\exists^Q_{F d_i}$ denotes the left adjoint to the map $\pi_{F d_i}^\ast \colon  Q( d_i \times c) \to Q(c)$, and so it remains to show that this yields a well-defined morphism of $\{\exists,\land,\lor\}$-doctrines.

First, it is easily observed that $\alpha_c^\exists$ preserves joins:
\begin{align*}
	& \alpha_c^\exists(\{({d_1},x_1) , \dots , (d_n,x_n)\} \lor \{(e_1,y_1) , \dots , (e_m,y_m)\}) \\
	= \ & \alpha_c^\exists(\{({d_1},x_1) , \dots , (d_n,x_n),(e_1,y_1) , \dots , (e_m,y_m)\}), \\
	= \ & \exists^Q_{Fd_1} \alpha_{d_1 \times c}(x_1) \lor \dots \lor \exists^Q_{Fd_n} \alpha_{d_n \times c}(x_n) \lor \exists_{Fe_1}^Q \alpha_{e_1 \times c}(y_1) \lor \dots \lor \exists_{Fe_m}^Q \alpha_{e_m \times c}(y_m), \\
	= \ & \alpha_c^\exists(\{({d_1},x_1) , \dots , (d_n,x_n)\}) \lor \alpha_c^\exists(\{(e_1,y_1) , \dots , (e_m,y_m)\}).
\end{align*}
As before, meets are more involved but nonetheless straightforward:
\begin{align*}
	& \alpha_c^\exists(\set{(d_i,x_i)}{i \leqslant n} \land  \set{(e_j,y_j)}{j \leqslant m}) \\
	= \ & \alpha_c^\exists(\set{(d_i \land e_j , {\pi'}_{e_j}^\ast x_i \land {\pi'}_{d_i}^\ast y_j)}{i \leqslant n, j \leqslant m}), \\
	= \ & \textstyle \bigvee_{i \leqslant n, j \leqslant m} \exists_{F(d_i \times e_j)}^Q \alpha_{d_i \times e_j \times c}({\pi'}_{e_j}^\ast x_i \land {\pi'}_{d_i}^\ast y_j), \\
	= \ & \textstyle \bigvee_{i \leqslant n, j \leqslant m} \exists_{Fd_i \times Fe_j}^Q \alpha_{d_i \times e_j \times c}({\pi'}_{e_j}^\ast x_i \land {\pi'}_{d_i}^\ast y_j), \\
	= \ & \textstyle \bigvee_{i \leqslant n, j \leqslant m} \exists_{Fd_i }^Q \exists_{ Fe_j}^Q \alpha_{d_i \times e_j \times c}({\pi'}_{e_j}^\ast x_i \land {\pi'}_{d_i}^\ast y_j), \\
	= \ & \textstyle \bigvee_{i \leqslant n, j \leqslant m}  \exists_{Fd_i }^Q \exists_{ Fe_j}^Q ( {\pi'}_{e_j}^\ast \alpha_{d_i \times c}(x_i) \land {\pi'}_{d_i}^\ast \alpha_{e_j \times c}( y_j) ) & \text{(naturality of $\alpha$),} \\
	= \ & \textstyle \bigvee_{i \leqslant n, j \leqslant m}  \exists_{Fd_i }^Q ( \alpha_{d_i \times c}(x_i) \land \exists_{Fe_j}^Q {\pi'}_{d_i}^\ast \alpha_{e_j \times c}( y_j)  ) & \text{(Frobenius)}. \\
	= \ & \textstyle \bigvee_{i \leqslant n, j \leqslant m}  \exists_{Fd_i }^Q ( \alpha_{d_i \times c}(x_i) \land {\pi''}_{d_i}^\ast \exists^Q_{Fej} \alpha_{e_j \times c}(y_j) & \text{(Beck-Chevalley)}, \\
	= \ & \textstyle \bigvee_{i \leqslant n, j \leqslant m}  \exists_{Fd_i }^Q  \alpha_{d_i \times c}(x_i) \land \exists^Q_{Fej} \alpha_{e_j \times c}(y_j) & \text{(Frobenius)}, \\
	= \ & \textstyle \left( \bigvee_{i \leqslant n} \exists_{Fd_i}^Q \alpha_{d_i \times c} (x_i) \right) \land \left( \bigvee_{j \leqslant m} \exists_{Fe_j}^Q \alpha_{e_j \times c}(y_j)\right) , \\
	= \ & \alpha_c^\exists(\set{(d_i,x_i)}{i \leqslant n}) \land  \alpha_c^\exists( \set{(e_j,y_j)}{j \leqslant m}).
\end{align*}
Finally, that the square
\[
\begin{tikzcd}
	P^\exists(d \times c) \ar{r}{\Sigma_d} \ar{d}[']{\alpha_{d \times c}^\exists} & P^\exists(c) \ar{d}{\alpha^\exists_c} \\
	Q(F(d\times c)) \ar{r}{\exists_{Fd}^Q} & Q(Fc)
\end{tikzcd}
\]
commutes for all $d \in \cat$ is easily calculated:
\begin{align*}
	& \alpha^\exists_c ( \Sigma_d (\{(e_1,x_1), \dots , (e_n,x_n)\})) \\
	= \ & \alpha^\exists_c ( \{(d \times e_1, x_1), \dots , (d \times e_n, x_n)\} ), \\
	= \ & \exists_{F(d \times e_1)}^Q \alpha_{d \times e_1 \times c}(x_1) \lor \dots \lor \exists_{F(d \times e_n)}^Q \alpha_{d \times e_n \times c}(x_n), \\
	= \ & \exists_{Fd \times Fe_1}^Q \alpha_{d \times e_1 \times c}(x_1) \lor \dots \lor \exists_{Fd \times F e_n}^Q \alpha_{d \times e_n \times c}(x_n), \\
	= \ & \exists_{Fd}^Q \exists_{Fe_1}^Q \alpha_{d \times e_1 \times c}(x_1) \lor \dots \lor \exists_{Fd}^Q \exists_{F e_n}^Q \alpha_{d \times e_n \times c}(x_n), \\
	= \ & \exists_{Fd}^Q( \exists_{Fe_1}^Q \alpha_{d \times e_1 \times c}(x_1) \lor \dots \lor \exists_{F e_n}^Q \alpha_{d \times e_n \times c}(x_n) ), \\
	= \ & \exists_{Fd}^Q \alpha_{d \times c}^\exists ( \{ (e_1,x_1) , \dots , (e_n,x_n) \}),
\end{align*}
where we have used that $\exists_{Fd}^Q$, being a left adjoint, commutes with joins.
\endproof
This essentially completes the proof of Theorem \ref{thm:exist_comp_for_coherent}, in that Proposition \ref{prop:exist_comp_1cat} shows that the forgetful 1-functor $\Doc{\exists,\land,\lor} \hookrightarrow \Doc{\land,\lor}$ has a left 1-adjoint.  To deduce the 2-categorical version stated in Theorem \ref{thm:exist_comp_for_coherent}, we need only note that if $\omega \colon F \Rightarrow G$ is the underlying natural transformation of a transformation between a pair of morphisms of $\{\land,\lor\}$-doctrines $(F,\alpha), (G,\beta) \colon P \rightrightarrows Q$ where $Q$ is a $\{\exists,\land,\lor\}$-doctrine, then $\omega$ is also the underlying natural transformation of a transformation of the extended doctrine morphisms $\omega \colon (F,\alpha^\exists) \Rightarrow (G,\beta^\exists)$.
\begin{lem}
	Let $(F,\alpha) , (G,\beta)\colon P \rightrightarrows Q$ be a pair of morphisms of $\{\land,\lor\}$-doctrines and suppose that $Q$ is a $\{\exists,\land,\lor\}$-doctrine.  Given a natural transformation $\omega \colon F \Rightarrow G$, then $\alpha_c(x) \leqslant \omega_c^\ast \beta_c(x)$, for all $ c \in \cat$ and $x \in P(c)$, if and only if
	\[
	\alpha_c^\exists (\{(d_1,x_1) , \dots , (d_n,x_n)\}) \leqslant \omega_c^\ast \beta^\exists_c(\{(d_1,x_1), \dots , (d_n,x_n)\})
	\]
	for all $c \in \cat$ and $\{(d_1,x_1), \dots , (d_n,x_n)\} \in P^\exists(c)$.  
\end{lem}
\proof
Suppose that $\alpha_c(x) \leqslant \omega_c^\ast \beta_c(x)$, for all $ c \in \cat$ and $x \in P(c)$.  We first show that this implies that 
\[\alpha_c^\exists(\{(d,x)\}) = \exists^Q_{Fd} \alpha_{d \times c}(x) \leqslant \omega_c^\ast \exists_{Gd}^Q \beta_{d \times c}(x) =  \omega_c^\ast \beta_c^\exists(\{(d,x)\}).\]
First note that, by the fact that $F$ and $G$ preserve products and $\omega$ is natural, the component $\omega_{d \times c} \colon F d \times F c \to G d \times G c$ is the uniquely induced morphism making the diagram
\[
\begin{tikzcd}
	F d \ar{d}[']{\omega_d} & \ar{l} \ar{r} Fd \times F c \ar[dashed]{d}{\omega_{c \times d}} & F c \ar{d}{\omega_d} \\
	G d & \ar{l} \ar{r} G d \times G c & G(c)
\end{tikzcd}
\]
commute, and so $\omega_{d \times c} = (\omega_d \times 1_{G c}) \circ (1_{F d} \times \omega_c)$.  Therefore, we calculate that
\begin{align*}
	\alpha_{d \times c}(x) & \leqslant \omega_{d \times c}^\ast \beta_{d \times c}(x) = (1_{Fd} \times \omega_c)^\ast (\omega_d \times 1_{G c})^\ast \beta_{d \times c}(x), \\
	\implies \exists_{Fd}^Q \alpha_{d \times c}(x) & \leqslant \exists_{F d}^Q (1_{Fd} \times \omega_c)^\ast (\omega_d \times 1_{G c})^\ast \beta_{d \times c}(x), \\
	 & = \omega_c^\ast \exists_{F d}^Q  (\omega_d \times 1_{G c})^\ast \beta_{d \times c}(x) & \text{(Beck-Chevalley)} \\
	& \leqslant \omega_c^\ast \exists_{Gd}^Q \beta_{d \times c}(x) =  \omega_c^\ast \beta_c^\exists(\{(d,x)\}),
\end{align*}
where in the last deduction we have used the canonical inequality $\exists_{Fd} (\omega_d \times 1_{Gc})^\ast \leqslant 1_{Gc}^\ast \exists_{Gd}^Q = \exists_{Gd}^Q$, the \emph{mate}, coming from the commutativity of the square on the right, a consequence of the commutativity of the square on the left:
\[
\begin{tikzcd}
	F d \times G c \ar{d}[']{\pi_{Fd}} \ar{r}{ \omega_d \times 1_{Gc}}  & G d \times G c \ar{d}{\pi_{Gd}} \\
	Gc \ar[equal]{r} & G c,
\end{tikzcd}
\begin{tikzcd}
	Q( F d \times G c) & Q(G d \times Gc) \ar{l}[']{(\omega_d \times 1_{Gc})^\ast} \\
	Q(Gc) \ar{u}{\pi_{Fd}^\ast} & Q(Gc). \ar{l}{1_{Gc}^\ast} \ar{u}[']{\pi_{Gd}^\ast} 
\end{tikzcd}
\]
For an arbitrary $\{(d_1,x_1),\dots,(d_n,x_n)\} \in P^\exists$, the inequality
\begin{align*}
	\alpha_c^\exists(\{({d_1},x_1) , \dots , (d_n,x_n)\}) &  = \exists^Q_{Fd_1} \alpha_{d_1 \times c}(x_1) \lor \dots \lor \exists^Q_{Fd_n} \alpha_{d_n \times c}(x_n) \\
& \leqslant \omega_c^\ast  \beta^\exists_c(\{(d_1,x_1), \dots , (d_n,x_n)\}) \\
 & = \omega_c^\ast( \exists^Q_{Gd_1} \beta_{d_1 \times c}(x_1) \lor \dots \lor \exists^Q_{Gd_n} \beta_{d_n \times c}(x_n)), \\
 & = \omega_c^\ast  \exists^Q_{Gd_1} \beta_{d_1 \times c}(x_1) \lor \dots \lor \omega_c^\ast \exists^Q_{Gd_n} \beta_{d_n \times c}(x_n)
\end{align*}
is a consequence of the inequality $\exists_{Fd_i}^Q \alpha_{d_i \times c}(x_i) \leqslant \omega_c^\ast \exists_{Gd_i}^Q \beta_{d_i \times c}(x)$ established above.

For the converse implication, we would have, as a particular case, the desired inequality
\[ \alpha_c(x) = \alpha_c^\exists \eta_c(x) = \alpha^\exists(\{(1,x_1)\}) \leqslant \omega_c^\ast \beta^\exists_c(\{(1,x)\}) =  \omega_c^\ast \beta_c^\exists \eta_c(x) = \omega_c^\ast \beta_c(x).\]
\endproof
\begin{coro}\label{coro:statement_of_2adj}
	For each $\{\land,\lor\}$-doctrine $P$ and each $\{\exists,\land,\lor\}$-doctrine $Q$, there is a natural isomorphism of categories $\Doc{\land,\lor}(P,Q) \cong \Doc{\exists,\land,\lor}(P^\exists,Q)$.
\end{coro}
\subsection{The existential completion respects elementary structure}
Just as in \cite[\S 6]{trotta}, we show that the existential completion for $\{\land,\lor\}$-doctrines respects elementary structure when it is present, i.e.\ the existential completion studied above restricts to a left 2-adjoint to the forgetful functor $\Doc{\exists,\land,\lor,=} \hookrightarrow \Doc{\land,\lor,=}$.
\begin{prop}
	If $P$ is a $\{\land,\lor,=\}$-doctrine, then $P^\exists $ is a $\{\exists,\land,\lor,=\}$-doctrine, and moreover the equality predicate, for an object $d \in \cat$, is given by $\{(1,\delta_d)\}$.
\end{prop}
\proof
We need to show that ${\pi_d^\exists}^\ast(-) \land {\pi_{d'}^\exists}^\ast \{(1,\delta_d)\} \colon P^\exists(d' \times d) \to P^\exists(d' \times d \times d)$ is left adjoint to ${(1_{d'} \times \Delta_d)^\exists}^\ast \colon P^\exists(d' \times d \times d) \to P^\exists(d' \times d)$.  Given $\set{({e_i},x_i)}{i \leqslant n} \in P^\exists(d' \times d)$, we have that:
\begin{align*}
&	{\pi_d^\exists}^\ast(\set{({e_i},x_i)}{i \leqslant n}) \land {\pi_{d'}^\exists}^\ast \{(1,\delta_d)\} \\
= \  & \set{({e_i},(1_{e_i} \times \pi_d)^\ast x_i)}{i \leqslant n} \land \{(1,(1_1 \times \pi_{d'})^\ast \delta_d)\}, \\
= \  & \set{({e_i},(1_{e_i} \times \pi_d)^\ast x_i \land \pi_{d'}^\ast \delta_d )}{i \leqslant n} ,
\end{align*}
and so, for $\set{({b_p},y_p)}{p \leqslant q} \in P^\exists(d' \times d \times d)$,
\[{\pi_d^\exists}^\ast(\set{({e_i},x_i)}{i \leqslant n}) \land {\pi_{d'}^\exists}^\ast \{(1,\delta_d)\}  \leqslant  \set{({b_p},y_p)}{p \leqslant q},\]
if and only if, for each $i \leqslant n$, there exists a finite family of arrows
\[
r_1 \colon e_i \times d' \times d \times d \to b_{j_1} \times d' \times d \times d, \dots , r_k \colon e_i \times d' \times d \times d \to b_{j_k} \times d' \times d \times d
\]
satisfying the commutativity condition that $\pi'_{e_i} = \pi_{b_{j_\ell}} \circ r_\ell$, for each $\ell \leqslant k$, and such that $(1_{e_i} \times \pi_d)^\ast x_i \land \pi_{d'}^\ast \delta_d \leqslant r_1^\ast y_{j_1} \lor \dots \lor r_k^\ast y_{j_k}$.  Using that $\delta_d$ is the equality predicate in the elementary doctrine $P \colon \cat\op \to \DLat$, this latter inequality holds if and only if
\begin{align*}
	x_i \leqslant & (1_{e_i \times d} \times \Delta_d)^\ast ( r_1^\ast y_{j_1} \lor \dots \lor r_k^\ast y_{j_k} ), \\
	= & (1_{e_i \times d} \times \Delta_d)^\ast r_1^\ast y_{j_1} \lor \dots \lor (1_{e_i \times d} \times \Delta_d)^\ast r_k^\ast y_{j_k}.
\end{align*}

Meanwhile, we have that
\[
\set{({e_i},x_i)}{i \leqslant n} \leqslant {(1_{d'} \times \Delta_d)^\exists}^\ast \set{({b_p},y_p)}{p \leqslant q} = \set{({b_p},(1_{b_p} \times 1_{d'} \times \Delta_d)^\ast y_p)}{p \leqslant q}
\]
if and only if, for each $i \leqslant n$, there exists a finite family of arrows
\[
r'_1 \colon e_i \times d' \times d \to b_{j'_1} \times d' \times d, \dots , r'_{k'} \colon e_i \times d' \times d \to b_{j'_{k'}} \times d' \times d
\]
such that $\pi_{e_i} = \pi'_{b_p} \circ r'_\ell$, for each $\ell \leqslant k'$, and 
\[x_i \leqslant {r'_1}^\ast (1_{b_{j'_1} \times d'} \times \Delta_d)^\ast y_{j'_1} \lor \dots \lor {r'_1}^\ast (1_{b_{j'_{k'}} \times d'} \times \Delta_d)^\ast y_{j'_{k'}}.\]

By manipulating the universal property of the product, whenever we are given an arrow $r_\ell \colon e_i \times d' \times d \times d \to b_{j_\ell} \times d' \times d \times d$, we can easily obtain an arrow $r'_\ell \colon e_i \times d' \times d \to b_{j_\ell} \times d' \times d$ making the square commute
\[
\begin{tikzcd}
	e_i \times d' \times d \ar[dashed]{r}{r'_\ell} \ar{d}[']{1_{e_i \times d'} \times \Delta_d} & b_{j_\ell} \times d' \times d \ar{d}{1_{b_{j_\ell} \times d'} \times \Delta_{d}} \\
	e_i \times d' \times d \times d \ar{r}{r_\ell} & b_{j_\ell} \times d' \times d \times d,
\end{tikzcd}
\]
and, {\it vice versa}, starting with such an arrow $r'_\ell$ as above, we can obtain an arrow $r_\ell$ making the above square commute. Therefore, given, for each $i \leqslant n$,  a finite family of arrows $\set{r_\ell \colon e_i \times d' \times d \times d \to b_{j_\ell} \times d' \times d \times d}{\ell \leqslant k}$ witnessing that 
\[{\pi_d^\exists}^\ast(\set{({e_i},x_i)}{i \leqslant n}) \land {\pi_{d'}^\exists}^\ast \{(1,\delta_d)\}  \leqslant  \set{({b_p},y_p)}{p \leqslant q},\]
we can obtain a family of arrows $\set{r'_\ell \colon e_i \times d' \times d \to b_{j_\ell} \times d' \times d}{\ell \leqslant k}$ that satisfies the inequality
\begin{align*}
		x_i & \leqslant  (1_{e_i \times d} \times \Delta_d)^\ast r_1^\ast y_{j_1} \lor \dots \lor (1_{e_i \times d} \times \Delta_d)^\ast r_k^\ast y_{j_k}, \\
		& = {r'_1}^\ast (1_{b_{j_1} \times d'} \times \Delta_d)^\ast y_{j_1} \lor \dots \lor {r'_k}^\ast (1_{b_{j_{k}} \times d} \times \Delta_d)^\ast y_{j_{k}}.
\end{align*}
Hence, we have that
\begin{align*}
	{\pi_d^\exists}^\ast(\set{({e_i},x_i)}{i \leqslant n}) \land {\pi_{d'}^\exists}^\ast \{(1,\delta_d)\}  & \leqslant  \set{({b_p},y_p)}{p \leqslant q} \\
	\implies  \quad \set{({e_i}, x_i)}{i \leqslant n} & \leqslant {(1_{d'} \times \Delta_d)^\exists}^\ast \set{({b_p},y_p)}{p \leqslant q},
\end{align*}
and via a similar argument, {\it mutatis mutandis}, the converse implication also holds, proving the desired adjunction.
\endproof
Thus, since the equality predicate for the existential completion $P^\exists \colon \cat\op \to \DLat$, at an object $d \in \cat$, is given by $\eta_d(\delta_d) = \{(1,\delta_d)\}$, we immediately deduce the following, which, combined with Corollary \ref{coro:statement_of_2adj}, constitutes a proof of Proposition \ref{prop:univ_prop_coh_and_equality}.
\begin{coro} Whenever $P \colon \cat\op \to \DLat$ is a $\{\land,\lor,=\}$-doctrine, 
	\begin{enumerate}
		\item The unit $\eta \colon P \to P^\exists$ is a morphism of $\{\land,\lor,=\}$-doctrines,
		\item If $(F,\alpha) \colon P \to Q$ is a morphism of $\{\land,\lor,=\}$-doctrines, the unique extension 
		\[(F,\alpha^\exists) \colon P^\exists \to Q\]
		constructed in Proposition \ref{prop:exist_comp_1cat} is a morphism of $\{\exists,\land,\lor,=\}$-doctrines.
	\end{enumerate}
\end{coro}

%
%

\printbibliography

\end{document}